\documentclass[final,3p,sort&compress]{elsarticle}

\usepackage{hyperref}
\usepackage{graphicx}
\usepackage{xcolor}
\usepackage{amsmath}
\usepackage{mathtools}
\usepackage{amsthm}
\usepackage{yhmath}
\usepackage{mathrsfs}
\usepackage{delimset}
\usepackage{empheq}
\usepackage{units}
\usepackage{amssymb}
\usepackage{stmaryrd}

\everydisplay{\let\binom\mybinom}

\allowdisplaybreaks

\makeatletter
\def\NAT@spacechar{}
\makeatother

\usepackage{titlesec}
\titleformat*{\subsection}{\bfseries}

\journal{Xxxxx}

\begin{document}

\begin{frontmatter}

\title{Tempered Fractional Brownian Motion with Variable Index and Variable Tempering Parameter
}

\author[label1]{S.C. Lim\corref{cor1}}
\address[label1]{50 Holland Road,
\#02-01 Botanika,
Singapore 258853
}
\ead{sclim47@gmail.com}
\cortext[cor1]{corresponding author}

\author[label5]{Chai Hok Eab}
\address[label5]{
7/9 Rongmuang 5, Pathumwan, Bangkok 10330,Thailand
}
\ead{Chaihok.E@gmail.com}

\begin{abstract}
Generalizations of tempered fractional Brownian from single index to two indices and variable index or tempered multifractional Brownian motion are studied. 
Tempered fractional Brownian motion and tempered multifractional Brownian motion with variable tempering parameter are considered.
\end{abstract}

\begin{keyword}
Fractional Ornstein-Uhlenbeck process, 
reduced fractional Ornstein-Uhlenbeck process, 
tempered fractional Brownian motion, 
variable tempering parameter,
tempered multifractional Brownian motion
\end{keyword}

\end{frontmatter}

\section{Introduction}
\label{sec:Introduciton}\noindent
Tempered fractional Brownian motion (TFBM) has been introduced recently by modifying the moving average representation of 
fractional Brownian motion (FBM) with the inclusion of an exponential tempering factor to the power-law kernel \cite{MeerschaertSbzikar2013,MeerschaertSabzikar2014} 
TFBM has found applications in in many phenomena, from transient anomalous diffusion \cite{ChenWeihuaDeng2017} wind speed and 
geophysical flow \cite{MeerschaertSabzikarM.PhanikumarEtAl2014,ZhangXiao2017}.

TFBM can be treated from the view point of fractional Ornstein-Uhlenbeck process (FOU). 
Recall that in their definition of FBM, denoted by $B_H(t)$,
Mandelbrot and van Ness \cite{MandelbrotVanNess68} introduced the reduced form of FBM $B_H(t)=X_H(t) - X_H(0)$ 
to get rid of divergence in the original definition of FBM $X_H(t)$ based on Weyl fractional integration of white noise. 
Although the FOU is well-defined using Weyl fractional integral, 
its covariance diverges when the damping constant $\lambda$ goes to 0. 
The $\lambda \to 0$  limit of the reduced process of FOU is well-defined, 
and it can be shown \mbox{\cite{LimTeo2007}} to be the FBM of Mandelbrot and van Ness.

The main aim of this paper is to study some generalizations of TFBM in terms of the RFOU. 
The advantages of adopting such an approach are two folds. 
First, the known properties of FOU allow one to verify various properties of TFBM directly or with some modifications. 
Another advantage of treating TFBM as reduced process of FOU is that 
it allows extension of TFBM to two indices in a more direct way. 
Other generalisations, 
which include TFBM with time-dependent tempering, 
and variable index TFBM or tempered multifractional Brownian motion (TMBM), 
and TMBM with variable tempering parameter, 
can be carried out based on their corresponding reduced processes. 
To the best of authors' knowledge, 
these generalisations have not been studied before. 
Some of the nice properties of these processes allow more realistic and flexible modelling of various natural and man-made phenomena.

The plan of this paper is as follows. Section \ref{sec:RfOU} consider TFBM as the reduced FOU process. 
Extension of TFBM to two indices is studied in section \ref{sec:.TfOU2ind}.
The subsequent sections discuss the generalization of TFBM to tempered multifractional Brownian motion (TMBM) with variable index, 
and the extension of TFBM and TMBM with variable tempering parameer. 
Concluding remarks are given in the final section.

\section{Tempered Fractional Brownian Motion as Reduced Fractional Ornstein-Uhlenbeck Process}
\label{sec:RfOU}\noindent
Recall that FOU can be obtained as the solution of the following fractional Langevin equation:
\begin{align}
\bigl({}_aD_t + \lambda\bigr)^\alpha X_{\alpha,\lambda}(t) & = \eta(t), &
    \alpha & > 0,
\label{eq:RfOU_0010}
\end{align}
where $\lambda > 0$ is a positive constant, 
and $\eta(t)$ is the Gaussian white noise with zero mean and delta-correlated covariance. 
The fractional derivative ${}_aD_t$  is defined by \cite{Samko1993,Kilbas2006}
\begin{align}
  {_aD_t^\alpha}f(t) & = \frac{1}{\Gamma(n-\alpha)} \left(\frac{d}{dt}\right)^n  \int_a^t \frac{f(u)}{(t-u)^{\alpha-n+1}} du, &
          & n-1 < \alpha < n. 
\label{eq:RfOU_0020}
\end{align}
For $a=0 $, the fractional derivative is known as the Riemann-Liouville fractional derivative; 
when $a=-\infty$, ${_\infty D_t^\alpha}$  is called the Weyl fractional derivative. 

The following operator identity holds for both the Riemann-Liouville and Weyl fractional derivatives \cite{Eab06b}
\begin{align}
   \bigl({_aD_t} + \lambda\bigr)^\alpha & = e^{-\lambda{t}} {_aD_t^{\alpha}} e^{\lambda{t}} .
\label{eq:RfOU_0030}
\end{align}
The shifted fractional derivative (\ref{eq:RfOU_0030}) is also known as tempered fractional derivative in subsequent work 
\cite[see for examples, ][]{CarteaCastillo-Negrete2007,SabzikarMeerschaertChen2015,CaoLiChen2014}.

In this paper only FOU of Weyl type would be considered. 
One has for $\alpha > 1/2$,
\begin{subequations}
\label{eq:RfOU_0040}
\begin{align}
  X_{\alpha,\lambda}(t) & = \frac{1}{\Gamma(\alpha)}\int_{-\infty}^t e^{-\lambda(t-u)}(t-u)^{\alpha-1}\eta(t)du,  
\label{eq:RfOU_0040a} \\
                     & = \frac{1}{\Gamma(\alpha)}\int_{-\infty}^\infty e^{-\lambda(t-u)_{+}}(t-u)_{+}^{\alpha-1}\eta(t)du,
\label{eq:RfOU_0040b}
\end{align}  
\end{subequations}
where $(x)_{+}^\mu =\bigl(\max(x,0)\bigr)^\mu$ with $0^0 = 0$.
Note that the condition $\alpha > 1/2$ is to ensure the above integral exists. 
(\ref{eq:RfOU_0040}) is known as the moving average representation FOU of Weyl type. 

$X_{\alpha,\lambda}(t)$ is a centred Gaussian stationary process with the following covariance and variance
\begin{align}
  C_{\alpha,\lambda}(t-s) & = \Bigl<X_{\alpha,\lambda}(t)X_{\alpha,\lambda}(s)\Bigr> 
                         =  \frac{1}{\sqrt{\pi}\Gamma(\alpha)}\left(\frac{|t-s|}{2\lambda}\right)^{\alpha-1/2} K_{\alpha-1/2}\bigl(\lambda|t-s|\bigr) , 
\label{eq:RfOU_0050}\\
  \sigma_{\alpha,\lambda}^2(t) & = \Bigl<\bigl(X_{\alpha,\lambda}(t)\bigr)^2\Bigr> 
                             = \frac{\Gamma(2\alpha-1)}{\bigl(\Gamma(\alpha)\bigr)^2(2\lambda)^{2\alpha-1}} .
\label{eq:RfOU_0060}
\end{align}
$K_\nu(z)$ is the modified Bessel function of second kind \cite{AbramowitzStegun64}.
The spectral density of $X_{\alpha,\lambda}(t)$ is
\begin{align}
  S_{\alpha,\lambda}(k) & = \frac{1}{2\pi}\int_{-\infty}^\infty C_{\alpha,\lambda}(\tau)e^{ik\tau} d\tau = \frac{1}{2\pi}\frac{1}{\bigl(k^2 + \lambda^2\bigr)^\alpha}.
\label{eq:RfOU_0070}
\end{align}
FOU has the following spectral representation
\begin{align}
  X_{\alpha,\lambda}(t) & = \frac{1}{\sqrt{2\pi}}\int_{-\infty}^\infty \frac{e^{-ikt}\widetilde{\eta}(k)dk}{(-ik+\lambda)^\alpha}.
\label{eq:RfOU_0080}
\end{align}
$X_{\alpha,\lambda}(t)$ is a stationary process with short-range dependent \cite{Eab06b,Eab06a,LimLITeo2007}.

In contrast to FBM, (\ref{eq:RfOU_0040}) provides a well-defined FOU based on the Weyl-fractional integral. 
However, the covariance and variance of FOU of Weyl type are divergent in the $\lambda \to 0$ limit \cite{LimTeo2007}. 
The divergence disappears if one considers the reduced version of FOU, 
analogous to the reduced process introduced for FBM \cite{MeerschaertSbzikar2013}.

TFBM is defined as the reduced process of FOU for $\alpha > 1/2$:
\begin{align}
B_{\alpha,\lambda}(t) & = X_{\alpha,\lambda}(t) - X_{\alpha,\lambda}(0) \nonumber \\
                   & = \frac{1}{\Gamma(\alpha)}\int_{-\infty}^t e^{-\lambda(t-u)}(t-u)^{\alpha-1}\eta(t)du
                     - \frac{1}{\Gamma(\alpha)}\int_{-\infty}^0 e^{-\lambda(-u)}(-u)^{\alpha-1}\eta(t)du \nonumber \\
                   & = \frac{1}{\Gamma(\alpha)}\int_{-\infty}^\infty 
                                          \Bigl[
                                               e^{-\lambda(t-u)_{+}}(t-u)_{+}^{\alpha-1} - e^{-\lambda(-u)_{+}}(-u)_{+}^{\alpha-1}
                                          \Bigr]
                                          \eta(t)du.
\label{eq:RfOU_0090}
\end{align}

From (\ref{eq:RfOU_0070}) and (\ref{eq:RfOU_0080}) one gets the following spectral representation for TFBM:
\begin{align}
  B_{\alpha,\lambda}(t) & = \frac{1}{\sqrt{2\pi}}\int_{-\infty}^\infty \frac{\Bigl(e^{-ikt}-1\Bigr)\widetilde{\eta}(k)dk}{(-ik+\lambda)^\alpha}.
\label{eq:RfOU_0100}
\end{align}

\begin{figure}[ht]
  \centering
  \fbox{
  \includegraphics[bb = 0 0 350 280,width = 0.5\textwidth]{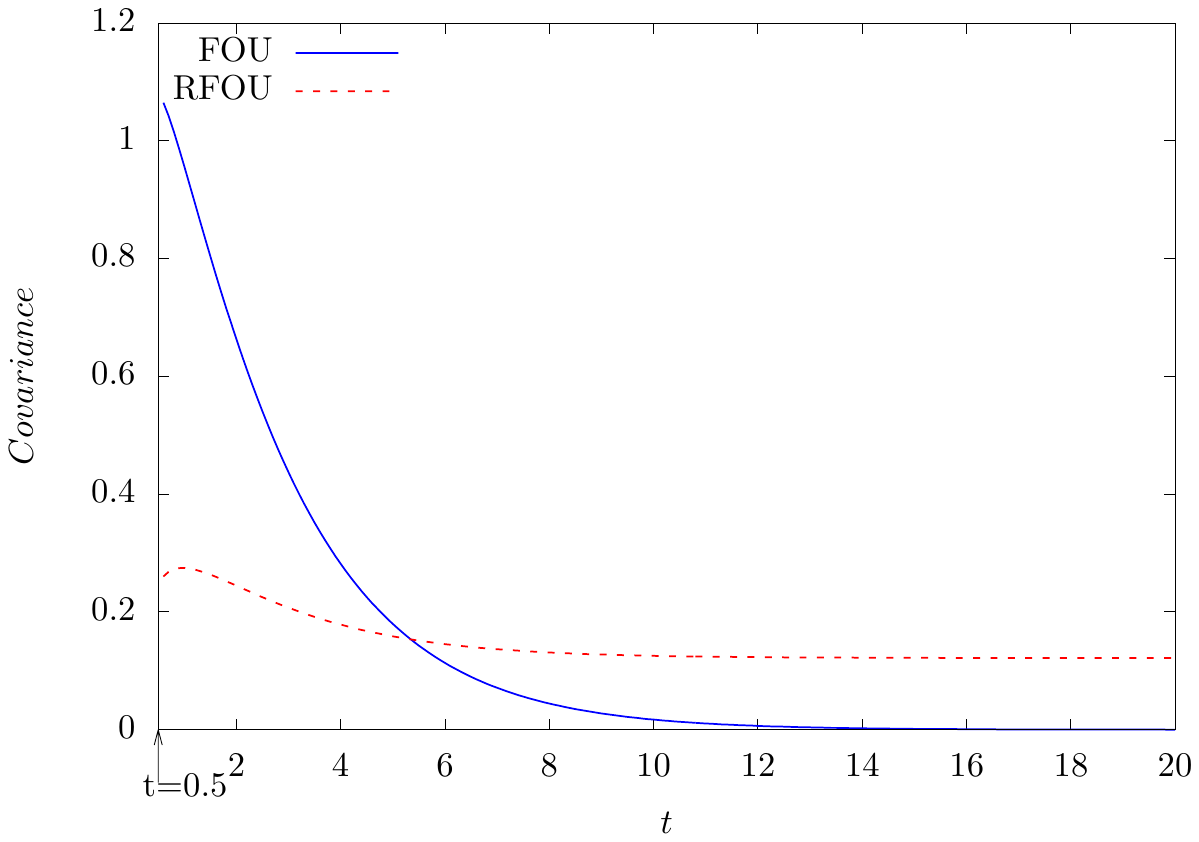}
  }
  \caption{covaraince of FOU and RFOU: $s=0.5$, $\lambda = 0.5, H = 0.75$}
  \label{fig:FOUvsTFBM_0010}
\end{figure}

The covariance of $B_{\alpha,\lambda}(t)$ can be calculated using (\ref{eq:RfOU_0050}) and (\ref{eq:RfOU_0060}):
\begin{align}
  \wideparen{C}_{\alpha,\lambda}(t,s) & = \bigl<B_{\alpha,\lambda}(t)B_{\alpha,\lambda}(s)\bigr> \nonumber \\
                               & = { \frac{1}{\sqrt{\pi}\Gamma(\alpha)}}
                                    \Biggl[
                                      \left(\frac{|t-s|}{2\lambda}\right)^{\alpha-1/2}
                                      K_{\alpha-1/2}(\lambda|t-s|)
                                    - \left(\frac{|t|}{2\lambda}\right)^{\alpha-1/2}
                                      K_{\alpha-1/2}(\lambda|t|)  \nonumber \\
                               & \qquad\qquad\qquad\qquad
                                    - \left(\frac{|s|}{2\lambda}\right)^{\alpha-1/2}
                                      K_{\alpha-1/2}(\lambda|s|)
                                    \Biggr]
                                 + \frac{\Gamma(2\alpha-1)}{\bigl(\Gamma(\alpha)\bigr)^2(2\lambda)^{2\alpha-1}}.
\label{eq:RfOU_0110}
\end{align}
The variance of $B_{\alpha,\lambda}(t)$ is
\begin{align} 
  \wideparen{\sigma}_{\alpha,\lambda}^2(t) & = \Bigl<\bigl(B_{\alpha,\lambda}(t)\bigr)^2\Bigr>
                                     =  \frac{2\Gamma(2\alpha-1)}{\bigl(\Gamma(\alpha)\bigr)^2(2\lambda)^{2H}}   
                                     -  2\frac{1}{\sqrt{\pi}\Gamma(\alpha)}\left(\frac{|t|}{2\lambda}\right)^{\alpha-1/2} K_{\alpha-1/2}(\lambda|t|).
\label{eq:RfOU_0120}
\end{align}
By letting $H=\alpha - 1/2$, 
\begin{align}
  c_t & = \frac{2\Gamma(2H)}{\bigl(\Gamma(H+1/2)\bigr)^2(2\lambda|t|)^{2H}} - \frac{2}{\sqrt{\pi}\Gamma(H+1/2)}\left(\frac{1}{2\lambda|t|}\right)^H K_H(\lambda|t|),
\label{eq:RfOU_0130}
\end{align}
the covariance of the TFBM can be expressed as
\begin{align}
  \wideparen{C}_{\alpha,\lambda}(t,s) & = \frac{1}{2}\Bigl[c_t|t|^{2H} + c_s|s|^{2H} - c_{t-s}|t-s|^{2H}\Bigr],
\label{eq:RfOU_0140}
\end{align}
which has the same form as the covariance of TFBM up to a multiplicative constant $\Gamma(H+1/2)^{-2}$ as given in \cite{MeerschaertSbzikar2013,MeerschaertSabzikar2014}.
Figure  \ref{fig:FOUvsTFBM_0010} shows the plots of covariance functions of FOU and its reduced process for 
$H = 0.75$ or $\alpha=1.25$, and $\lambda = 0.5$.

The following are the properties of TFBM can be obtained as direct consequences or some modifications of the properties of FOU. 

\subsection*{Stationarity property}
\label{sec:RfOU_Statn}\noindent
In contrast to FOU which is a stationary process, TFBM is non-stationary.
$B_{\alpha,\lambda}(t)$
is asymptotically stationary. 
The covariance of $B_{\alpha,\lambda}(t)$ can be expressed in the following form,
\begin{align}
  \wideparen{C}_{\alpha,\lambda}(t,t+\tau) & = C_{\alpha,\lambda}(\tau) - C_{\alpha,\lambda}(t+\tau) - C_{\alpha,\lambda}(t) + C_{\alpha,\lambda}(0) ,
\label{eq:RfOU_Statn_0010}
\end{align}
where $C_{\alpha,\lambda}(t)$ is the covariance of the FOU.
As $\tau \to \infty$, 
$C_{\alpha,\lambda}(\tau) - C_{\alpha,\lambda}(t+\tau) \to 0$ so $\wideparen{C}_{\alpha,\lambda}(t,t+\tau) \sim C_{\alpha,\lambda}(0) - C_{\alpha,\lambda}(t)$.
Thus, TFBM becomes a stationary process in the large time limit.

The increment $\Delta_\tau B_{\alpha,\lambda}(t) = B_{\alpha,\lambda}(t+\tau) - B_{\alpha,\lambda}(t) = X_{\alpha,\lambda}(t+\tau) - X_{\alpha,\lambda}(t)$ 
is a stationary process. 
This follows from the stationarity of $X_{\alpha,\lambda}(t)$.
The spectral representation of the increment process is
\begin{align}
  \Delta_\tau B_{\alpha,\lambda}(t) & = \frac{1}{\sqrt{2\pi}} 
                                      \int_{-\infty}^\infty 
                                      \frac{\Bigl(e^{-ik(t+\tau)} - e^{-ikt}\Bigr)\widetilde{\eta}(k)dk}
                                      {(-ik + \lambda)^\alpha} ,
\label{eq:RfOU_Statn_0020}
\end{align}
with covariance 
\begin{align}
  \wideparen{C}_{\Delta_\tau B_{\alpha,\lambda}}(t,s;\tau) 
          & = \frac{1}{2\pi} 
              \int_{-\infty}^\infty 
              \frac{\Bigl(e^{-ik(t+\tau)} - e^{-ikt}\Bigr)\Bigl(e^{ik(s+\tau)} - e^{iks}\Bigr)\widetilde{\eta}(k)dk}
              {\bigl(k^2 + \lambda^2\bigr)^\alpha}  \nonumber \\
          & = \frac{1}{\sqrt{\pi}\Gamma(\alpha)}
              \Biggl[
                2\left(\frac{|t-s|}{2\lambda}\right)^{\alpha-1/2}K_{\alpha-1/2}\bigl(\lambda|t-s|\bigr)
                -\left(\frac{|t-s+\tau|}{2\lambda}\right)^{\alpha-1/2}K_{\alpha-1/2}\bigl(\lambda|t-s+\tau|\bigr)
            \nonumber \\
          & \qquad\qquad
                -\left(\frac{|t-s-\tau|}{2\lambda}\right)^{\alpha-1/2}K_{\alpha-1/2}\bigl(\lambda|t-s-\tau|\bigr)
              \Biggr] .
\label{eq:RfOU_Statn_0030}
\end{align}

\subsection*{Scaling property}
\label{sec:RfOU_Scaling}
In contrast to FBM, $B_{\alpha,\lambda}(t)$ is not self-similar. 
However, it satisfies the following modified global scaling property:
\begin{align}
  B_{\alpha,\lambda}(rt) & \hat{=} r^{\alpha=1/2}B_{\alpha,r\lambda}(t),
\label{eq:RfOU_Scaling_0010}
\end{align}
where $r$ is a positive constant, 
and $B_{\alpha,r\lambda}(rt)$ is the same process as $B_{\alpha,\lambda}(rt)$ with $\lambda$ replaced by $r\lambda$.
Here $\hat{=}$ denotes equality in the sense of finite dimensional distributions.  
(\ref{eq:RfOU_Scaling_0010}) can be easily verified by showing
\begin{align}
  \bigl<B_{\alpha,\lambda}(rt)B_{\alpha,\lambda}(rs)\bigr> = r^{2\alpha-1} \bigl<B_{\alpha,r\lambda}(t)B_{\alpha,r\lambda}(s)\bigr>.
\label{eq:RfOU_Scaling_0020}
\end{align}
Note that the additional scaling is due to the presence of the tempering term $\lambda$, 
self-similarity is recovered when $\lambda=0$.

The above scaling property can loosely be regarded as a generalization of self-similarity property. 
One can easily show that such a scaling property also holds for stationary process FOU. 
Some examples of stationary processes which satisfy this scaling property include processes with the stretched exponential covariance
\cite{LimMunuandy2003}
\begin{align}
  C_{\nu,\lambda}^{SE}(t) & = \lambda^{-\nu} e^{-(\lambda|t|)^\nu}, &
  \nu & \in (0,2], &
   \lambda & > 0.
\label{eq:RfOU_Scaling_0030}
\end{align}
and processes with the generalized Cauchy-type covariance
\cite{LimTeo2009d}
\begin{align}
  C_{\beta,\lambda}^C(t) & = \lambda^{-\beta}\bigl(1+ \lambda^\beta |t|^\beta\bigr)^{-1}&
  \beta & \in (0,2], &
   \lambda & > 0.
\label{eq:RfOU_Scaling_0040}
\end{align}

Note that the non-stationary reduced processes associated to these processes also obey the same scaling behaviour. 
Thus, one may say that the scaling property (\ref{eq:RfOU_Scaling_0010}) is not unique to processes of Ornstein-Uhlenbeck type, 
which include FOU and RFOU or TFBM. 
It is satisfied by a wider class of stationary and non-stationary Gaussian processes. 
One common property for these processes is that the $\lambda \to 0$ limit of the covariance diverges, 
and for the reduced process its covariance becomes the covariance of FBM 
when $\lambda \to 0$ (up to a multiplicative constant). 
A process having such scaling property also satisfies locally self-similar property of order $\kappa = \alpha - 1/2$,
and its tangent process at a point $t_\circ$ given by FBM $B_H(t_\circ)$ with Hurst index $H = \alpha - 1/2$ in the case of TFBM \cite{LimEab2019}.
In other words, such a process behaves locally like FBM.

\subsection*{Fractal dimension}
\label{sec:RfOU_FD}
\noindent
Since the variance of the increment process of $B_{\alpha,\lambda}(t)$ satisfies
\begin{align}
  \wideparen{\sigma}^2_{\alpha,\Delta\tau}(t) & = \Bigl<\bigl(B_{\alpha,\lambda}(t+\tau) - B_{\alpha,\lambda}(t)\bigr)^2\Bigr>
                                                     \leq A|\tau|^{2\alpha-1},
\end{align}
thus, almost surely the sample path of $B_{\alpha,\lambda}(t)$ is H\"{o}lderian of order $(\alpha-1/2) -\epsilon$ for all $\alpha > 0$.
A process which is locally self-similar of order $\kappa > 0$ and its sample paths are a.s.\ $\kappa-\epsilon$.
H\"{o}lderian for all $\epsilon > 0$, then the fractal dimension of its graph is a.s. equals to $2-\kappa$
\cite{KentWood1997,Adler1981,BenassiCohenIstas2003,Falconer2003b}
Applying this result to $B_{\alpha,\lambda}(t)$ gives the fractal dimension $5/2 -\alpha$ or $2-H$ if $\alpha=H+1/2$.
Thus, both the TFBM and FBM have the same fractal dimension.

\subsection*{Long range dependence}
\label{sec:RfOU_LR}
\noindent
In contrast to FOU which is a short memory process, 
its reduced process or TFBM is long-range dependent. 
Consider the correlation function of a non-stationary Gaussian process $Z(t)$ with correlation function
\begin{align}
  R(t+\tau,t) & = \frac{C(t+\tau,t)}{\sqrt{\sigma^2(t+\tau)\sigma^2(t)}},
\label{eq:RfOU_LR_0010}
\end{align}
where $C(t+\tau,t)$ and $\sigma^2(t)$ are respectively the covariance and variance of $Z(t)$. 
The process is long range dependent if 
\begin{align}
  \int_0^\infty \bigl|R(t,t+\tau)\bigr| d\tau & = \infty,
\label{eq:RfOU_LR_0020}
\end{align}
otherwise it is short range dependent 
\cite{Beran1994,AyacheCohenLevyVehel2000,FandrinPAmblard2003}.
Alternatively, $Z(t)$ is said to have long range dependence if 
$R(t+\tau,t) \sim \tau^\gamma$, $- < \gamma < 0$ as $\tau \to \infty$, for all $t>0$.

From (\ref{eq:RfOU_Statn_0010}) one has for 
one has for $\tau \to \infty$, 
$\wideparen{C}_{\alpha,\lambda}(t+\tau,t) \sim C_{\alpha,\lambda}(0) - C_{\alpha,\lambda}(t) = \wideparen{\sigma}_{\alpha,\lambda}(t)/2$,
and
$K_\nu(\tau) \sim \sqrt{\pi/2\tau}e^{-\tau} \to 0$,
one get $\wideparen{\sigma}_{\alpha,\lambda}^2(t+\tau) \to \sigma_{\alpha,\lambda}^2(t)$.
Thus
\begin{align}
  \wideparen{R}_{\alpha,\lambda}(t+\tau,t) & = \frac{\wideparen{C}_{\alpha,\lambda}(t+\tau,t)}
                                             {\bigl[\wideparen{\sigma}_{\alpha,\lambda}^2(t)\wideparen{\sigma}_{\alpha,\lambda}^2(t+\tau)\bigr]^{1/2}} 
                                         \sim     
                                             \frac{1}{2}\sqrt{\frac{\wideparen{\sigma}_{\alpha,\lambda}^2(t)}{\sigma_{\alpha,\lambda}^2(t)}}
                                          > 0, &
                                        \tau & \to \infty,     
\label{eq:RfOU_LR_0030}   
\end{align}
which implies TFBM is a long memory process.

\section{Tempered Fractional Brownian Motion with Two Indices}
\label{sec:.TfOU2ind}
\noindent
The fractal dimension of the graph of the sample path of a stochastic process is a measure of its roughness. 
Long memory dependence of a stochastic process is associated with a heavy tail behaviour of the covariance function, 
which is also known as Hurst effect. 
For TFBM, just like FBM, both these two properties are characterised by Hurst index $H$. 
From theoretical point of view, fractal dimension and Hurst effect are independent of each other as fractal dimension is a local property, 
whereas long-memory dependence is a global characteristic. 
It is therefore useful to decouple the characterisation of these properties.

There are not many models of stochastic processes with covariance parametrized by two indices 
that allow independent estimation of fractal dimension and long memory dependence associated with the model. 
The few that have been studied include FOU with two indices \cite{LimLiTeo2008,LimTeo2009c} generalized Whittle-Matern \cite{LimTeo2009}
and generalized Cauchy models \cite{LimTeo2009d}.
In this section one more example, namely TFBM with two indices, is introduced based on the reduced process of FOU indexed by two parameters.

Fractional Langevin equation (\ref{eq:RfOU_0010}) can be extended to two indices 
\cite{LimLiTeo2008,LimTeo2009c}
\begin{align}
  \bigl({}_{-\infty}D_t^\beta + \lambda^\beta\bigr)X_{\alpha\beta,\lambda}(t) & = \eta(t), &
     \alpha\beta & > 1/2.
\label{eq:TfOU2ind_0010}
\end{align}
Here, $\lambda$ is replace by $\lambda^\beta$ to preserve the scaling property (\ref{eq:RfOU_Scaling_0010}).
By using Fourier transform method, the solution of (\ref{eq:TfOU2ind_0010}) is found to be a stationary Gaussian process
\begin{align}
  X_{\alpha\beta,\lambda}(t) & = \frac{1}{\sqrt{2\pi}} \int_{\mathbb{R}} \frac{e^{ikt}\widetilde{\eta}}{\bigl((ik)^\beta + \lambda^\beta\bigr)^\alpha} .
\label{eq:TfOU2ind_0020}
\end{align}
Again, the condition $\alpha\beta > 1/2$ is necessary to ensure the above integral is finite. 
$X_{\alpha\beta,\lambda}(t)$ has a more complicated spectral density,
\begin{align}
  S_{\alpha\beta,\lambda}(k) & = \frac{1}{2\pi\Bigl|(ik)^\beta + \lambda^\beta\bigr|^{2\alpha}}
                           = \frac{1}{2\pi\Bigl||k|^{2\beta} + 2\lambda^\beta|k|^\beta \cos(\alpha\pi/2) + \lambda^{2\beta}\bigr)^\alpha},
\label{eq:TfOU2ind_0030}
\end{align} 
which in general cannot be evaluated to give a closed analytic expression. 
Despite this, many of the basic properties of $X_{\alpha\beta,\lambda}(t)$ can still be obtained and studied \cite{LimLiTeo2008}.
Hence, one can consider the reduced process associated with it
\begin{align}
  B_{\alpha\beta,\lambda}(t) & = X_{\alpha\beta,\lambda}(t) - X_{\alpha\beta,\lambda}(0) 
                           = \frac{1}{\sqrt{2\pi}} \int_{-\infty}^\infty \frac{\bigl(e^{-ikt} -1\bigr)\widetilde{\nu}dk}{\bigl((-ik)^\beta + \lambda^\beta\bigr)^\alpha},
\label{eq:TfOU2ind_0040}
\end{align}
and examine its properties just like the case of single index. 
However, instead of (\ref{eq:TfOU2ind_0040}),
a different RFOU with two indices will be considered.

FOU with a simpler spectral density is given by the solution of the fractional Langevin equation of Riesz type with two indices 
\begin{align}
  \bigl(\mathbf{D}_t^{2\beta} + {\lambda}^{2\beta}\bigr)^{\alpha} Y_{\alpha\beta,\lambda}(t) & = \eta(t),
\label{eq:TfOU2ind_0050}
\end{align}
where $\mathbf{D}_t^\alpha$, $\alpha > 0$ is the one-dimensional Riesz derivative defined by 
\cite{Samko1993,Kilbas2006}
\begin{align}
  \mathbf{D}_t^\alpha f & = \left(-\frac{d^2}{dt^2}\right)^{\alpha/2}f 
                        = \mathcal{F}^{-1}\bigl(|k|^\alpha \widetilde{f}(k)\bigr),
\label{eq:TfOU2ind_0060}
\end{align}
where $\widetilde{f}(k)$ is the Fourier transform of $f(t)$.
For $\alpha\beta > 1/2$, the solution of (\ref{eq:TfOU2ind_0050}) is given by
\begin{align}
  Y_{\alpha\beta,\lambda}(t) & = \frac{1}{\sqrt{2\pi}}\int_{\mathbb{R}} \frac{e^{ikt} \widetilde{\eta}(k)}{\bigl(|k|^{2\beta}+\lambda^{2\beta}\bigr)^{\alpha/2}} dk.
\label{eq:TfOU2ind_0070}
\end{align}
The spectral density of the process has a simpler form as compared with (\ref{eq:TfOU2ind_0030}):
\begin{align}
  S_{\alpha\beta,\lambda}(k) & = \frac{1}{2\pi} \frac{1}{\bigl(|k|^{2\beta}+\lambda^{2\beta}\bigr)^\alpha}.
\label{eq:TfOU2ind_0080}
\end{align}
The covariance function $C_{\alpha\beta,\lambda}(t)$ can be obtained by taking the inverse Fourier transform of $S_{\alpha\beta,\lambda}(k)$.
However, it does not in general has a closed analytic form. 
The variance is
\begin{align}
  C_{\alpha\beta,\lambda}(0) & = \frac{1}{\pi}\int_0^\infty \frac{1}{\bigl(|k|^{2\beta}+\lambda^{2\beta}\bigr)^\alpha} dk
                           = \frac{\Gamma(1/2\beta)\Gamma(\alpha -1/2\beta)}{2\pi\beta\Gamma(\alpha)}\lambda^{1-2\alpha\beta}.
\label{eq:TfOU2ind_0090}
\end{align}
Note that both $X_{\alpha\beta,\lambda}(t)$  and  $Y_{\alpha\beta,\lambda}(t)$ 
can be regarded as two different generalisations of FOU of single index to two indices. 
These two processes have similar long and short time asymptotic properties \cite{LimLiTeo2008,LimTeo2009c}.
However, only the reduced process associated with  $Y_{\alpha\beta,\lambda}(t)$ will be considered here. 

Let $B_{\alpha\beta,\lambda}(t)$  denotes the reduced process associated with $Y_{\alpha\beta,\lambda}(t)$:
\begin{align}
  B_{\alpha\beta,\lambda}(t) & = Y_{\alpha\beta,\lambda}(t) - Y_{\alpha\beta,\lambda}(0). 
\label{eq:TfOU2ind_0100}
\end{align}
The spectral representation for $B_{\alpha\beta,\lambda}(t)$ is
\begin{align}
  B_{\alpha\beta,\lambda}(t) & = \frac{1}{\sqrt{2\pi}}\int_{\mathbb{R}} \frac{\bigl(e^{ikt} -1\bigr) \widetilde{\eta}(k)}{\bigl(|k|^{2\beta}+\lambda^{2\beta}\bigr)^{\alpha/2}} dk,
\label{eq:TfOU2ind_0110}
\end{align}
Its covariance and variance are respectively
\begin{align}
  \wideparen{C}_{\alpha\beta,\lambda}(t+\tau,t) & = \bigl<B_{\alpha\beta,\lambda}(t+\tau)B_{\alpha\beta,\lambda}(t)\bigr>
 = \frac{1}{2\pi}\int_{\mathbb{R}} \frac{\bigl(e^{ik|\tau|} - e^{ik|t+\tau|} -e^{ik|t|} +1\bigr)}{\bigl(|k|^{2\beta}+\lambda^{2\beta}\bigr)^\alpha} dk,
\label{eq:TfOU2ind_0120}
\end{align}
and
\begin{align}
  \wideparen{\sigma}_{\alpha\beta,\lambda}^2(t) & =  \frac{1}{\pi}\int_{\mathbb{R}} \frac{1-e^{ik|t|}}{\bigl(|k|^{2\beta}+\lambda^{2\beta}\bigr)^\alpha} dk,
                                             = 2\bigl(C_{\alpha\beta,\lambda}(0)-C_{\alpha\beta,\lambda}(t)\bigr) \nonumber \\
  & = \frac{\Gamma(1/2\beta)\Gamma(\alpha -1/2\beta)}{2\pi\beta\Gamma(\alpha)}\lambda^{1-2\alpha\beta}
      + \frac{1}{\pi}\int_0^\infty \frac{e^{ik|t|}}{\bigl(|k|^{2\beta}+\lambda^{2\beta}\bigr)^\alpha} dk .
\label{eq:TfOU2ind_0130}
\end{align}

The covariance of $B_{\alpha\beta,\lambda}(t)$ can again be expressed in the same form (\ref{eq:RfOU_0140}) for TFBM with single index, 
namely
\begin{align}
    \wideparen{C}_{\alpha\beta,\lambda}(t,s) & = \frac{1}{2}\Bigl[c_t|t|^{2H} + c_s|s|^{2H} - c_{t-s}|t-s|^{2H}\Bigr],
\label{eq:TfOU2ind_0140}
\end{align}
with
\begin{align}
  c_t & = \frac{\Gamma(1/2\beta)\Gamma(\alpha -1/2\beta)}{2\pi\beta\Gamma(\alpha)(\lambda|t|)^{2H}}
       - \frac{1}{2\pi|t|^{2H}}\int_{\mathbb{R}} \frac{e^{ik|t|}}{\bigl(|k|^{2\beta}+\lambda^{2\beta}\bigr)^\alpha} dk ,
\label{eq:TfOU2ind_0150}
\end{align}
and $H = \alpha\beta-1/2$.
Note that the second term in the coefficient $c_t$ does not have a closed form, one can obtain its asymptotic properties.

For $\beta \in (0,1)$, one has for $t \to \infty$,
\begin{align}
  \int_{\mathbb{R}}\frac{e^{ik|\tau|}}{\bigl(k^{2\beta}+\lambda^{2\beta}\bigr)^\alpha} dk
              & = 2Im \Bigg[
                      \frac{1}{\tau}
                      \int_0^\infty \frac{e^{-u}du}{\bigl(e^{i\beta\tau} (u/\tau)^{2\beta} + \lambda^{2\beta}\bigr)^\alpha}
                      \Biggr] \nonumber \\
              & = 2Im \Bigg[
                      \frac{1}{\tau}
                      \int_0^\infty e^{-u} \sum_{j=0}^\infty \frac{(-1)^j\Gamma(\alpha+j)}{j!\Gamma(\alpha)}
                         e^{-i\beta{j}\pi} \lambda^{-2\beta(\alpha+j)} \biggl(\frac{u}{\tau}\biggr)^{2\beta{j}}
                         du
                      \Biggr] \nonumber \\
              & \sim 2 \sum_{j=1}^\infty \frac{(-1)^j\Gamma(\alpha+j)\Gamma(2\beta{j}+1)}{j!\Gamma(\alpha)}
                         \sin(\beta{j}\pi) \lambda^{-2\beta(\alpha+j)}\tau^{-(2\beta{j}+1)}, 
\label{eq:TfOU2ind_0160}
\end{align}
with the leading term 
$2\alpha\Gamma(1+2\beta)\sin(\beta\pi)\lambda^{-2\beta(\alpha+1)} \tau^{-(2\beta+1)}$.
The long-time asymptotic depends only on $\beta$.
On the other hand, the small time asymptotic behaviour of $Y_{\alpha\beta,\lambda}(t)$, hence $B_{\alpha\beta,\lambda}(t)$,
depends on the arithmetic nature of both the parameters $\alpha$ and $\beta$ with rather complicated conditions \cite{KOH1995ab,ErdoganOstrovskii1998,LimTeo2010}.

For most practical purposes, 
it is sufficient to consider small time asymptotic behaviour of the covariance function for values of $\alpha$  and $\beta$  confined to $1/2 < \alpha\beta < 3/2$, 
or $0 < H < 1$ if $\alpha\beta = H + 1/2$.
For the short-time asymptotic behaviour consider first the variance of the increment process of $Y_{\alpha\beta,\lambda}(t)$.
One gets for $t \to 0$,
\begin{align}
  \wideparen{\sigma}_{\alpha\beta,\lambda}^2(t) & = \frac{1}{\pi}\int_{\mathbb{R}} \frac{1-\cos(k|t|)}{\bigl(|k|^{2\beta}+\lambda^{2\beta}\bigr)^\alpha} dk
                                                  \nonumber \\
     & = \frac{4}{\pi}\int_{\mathbb{R}} \frac{\sin^2(k|t|/2)}{\bigl(|k|^{2\beta}+\lambda^{2\beta}\bigr)^\alpha} dk
       = \frac{4|t|^{2\alpha\beta-1}}{\pi}\int_0^\infty \frac{\sin^2(k/2)}{\bigl(|k|^{2\beta}+\lambda^{2\beta}|t|^{2\beta}\bigr)^\alpha} dk \nonumber \\
     &  = \frac{4|t|^{2\alpha\beta-1}}{\pi}\int_0^\infty k^{-2\alpha\beta}\sin^2(k/2) dk + o\bigl(|t|^{2\alpha\beta -1}\bigr) \nonumber \\
     &  =  - \frac{|t|^{2\alpha\beta-1}}{\Gamma(\alpha\beta)\cos(\alpha\beta\pi)}+ o\bigl(|t|^{2\alpha\beta -1}\bigr), &
\text{as} \ t & \to 0.
\label{eq:TfOU2ind_0170}
\end{align}
Using (\ref{eq:TfOU2ind_0130}),
one gets for $t \to 0$,
\begin{align}
  \int_{\mathbb{R}} \frac{e^{ik|t|}}{\bigl(|k|^{2\beta}+\lambda^{2\beta}\bigr)^\alpha} dk 
  & = \frac{\Gamma(1/2\beta)\Gamma(\alpha -1/2\beta)}{\beta\Gamma(\alpha)}\lambda^{1-2\alpha\beta}
                 + \frac{\pi|t|^{2\alpha\beta-1}}{2\Gamma(\alpha\beta)\cos(\alpha\beta\pi)}+ o\bigl(|t|^{2\alpha\beta -1}\bigr).
\label{eq:TfOU2ind_0180}
\end{align}
Note that the small-time asymptotic behaviour depends on $\alpha$ and $\beta$ through the product $\alpha\beta$.
By replacing $\alpha$ by $\alpha/\beta$, the variance of the increment process of $B_{\alpha\beta,\lambda}(t)$ varies as $t^{-2\alpha-1}$ as $t \to 0$. 
On the other hand, the long-time asymptotic behaviour of the covariance varies as $t^{-(1+2\beta)}$,
which is independent of $\alpha$.
In contrast to FBM and TFBM, 
it is possible to separately characterize the short-time property such as fractal dimension, 
and the long-time behaviour like long-range dependence of $B_{\alpha\beta,\lambda}(t)$ by using two different parameters. 
$B_{\alpha\beta,\lambda}(t)$ satisfies the same properties as TFBM $B_{\alpha,\lambda}(t)$ with single index. 
First, we verify the scaling property for $Y_{\alpha\beta,\lambda}(t)$.
Consider the covariance of $Y_{\alpha\beta,\lambda}(rt)$, $r > 0$,
\begin{align}
 C_{\alpha\beta,\lambda}(r\tau) & =  \bigl<Y_{\alpha\beta,\lambda}(r(t+\tau))Y_{\alpha\beta,\lambda}(rt)\bigr> \nonumber \\
  & = \frac{1}{\pi}\int_0^\infty \frac{\cos(k|r\tau|)}{\bigl(|k|^{2\beta}+\lambda^{2\beta}\bigr)^\alpha} dk 
    = \frac{|r\tau|^{2\alpha\beta-1}}{\pi}\int_0^\infty  \frac{\cos(k)}{\bigl(k^{2\beta}+(\lambda|r\tau|)^{2\beta}\bigr)^\alpha} dk \nonumber \\
  & = \frac{r^{2\alpha\beta-1}|\tau|^{2\alpha\beta-1}}{\pi}\int_0^\infty \frac{\cos(k)}{\bigl(k^{2\beta}+(\lambda|r\tau|)^{2\beta}\bigr)^\alpha} dk 
    = r^{2\alpha\beta-1}\bigl<Y_{\alpha\beta,r\lambda}(t+\tau)Y_{\alpha\beta,r\lambda}(t)\bigr> ,
\label{eq:TfOU2ind_0190}
\end{align}
where $Y_{\alpha\beta,r\lambda}(t)$ is the same process as  $Y_{\alpha\beta,\lambda}(rt)$ with $\lambda$ replaced by $r\lambda$.
Using (\ref{eq:TfOU2ind_0120}) and (\ref{eq:TfOU2ind_0190}),
\begin{align}
  \wideparen{C}_{\alpha\beta,\lambda}(t,t+\tau) & = \bigl<B_{\alpha\beta,\lambda}(t)B_{\alpha\beta,\lambda}(t+\tau)\bigr>
                                              =  C_{\alpha\beta,\lambda}(\tau) {\color{red} -} C_{\alpha\beta,\lambda}(t+\tau) 
                                                  {\color{red} -} C_{\alpha\beta,\lambda}(t)              
                                                  + C_{\alpha\beta,\lambda}(0) ,          \nonumber    
\end{align}
one obtains the scaling property for $B_{\alpha\beta,\lambda}(t)$:
\begin{align}
  B_{\alpha\beta,\lambda}\bigl(r(t)\bigr) & = r^{\alpha\beta-1/2}B_{\alpha\beta,\lambda}(t).
\label{eq:TfOU2ind_0200}
\end{align}

By using the long-time asymptotic property of the covariance of $B_{\alpha\beta,\lambda}(t)$,
one can show that the process is asymptotically stationary. 
Its increment process is stationary, which follows from the stationary property of $Y_{\alpha\beta,\lambda}(t)$.

$B_{\alpha\beta,\lambda}(t)$ is a long memory process just like $B_{\alpha,\lambda}(t)$.
One can use the same argument as for $B_{\alpha,\lambda}(t)$.
For $\tau \to \infty$.
\begin{align}
  \wideparen{R}_{\alpha\beta,\lambda}(t,t+\tau) & \sim \frac{\wideparen{\alpha}_{\alpha\beta,\lambda}^2(t)}
                                           {2\sqrt{\wideparen{\sigma}_{\alpha\beta,\lambda}^2(t)\wideparen{\sigma}_{\alpha\beta,\lambda}^2(t+\tau)}}
                                     = \frac{1}{2} \sqrt{\frac{\wideparen{\sigma}_{\alpha\beta,\lambda}(t)}{\sigma_{\alpha\beta,\lambda}(t)}} ,                                           
\end{align}
By noting that one has $ \wideparen{R}_{\alpha\beta,\lambda}(t,t+\tau) > 0$,
hence $ \int_0^\infty \wideparen{R}_{\alpha\beta,\lambda}(t,t+\tau) d\tau = \infty$.
Thus, the process $B_{\alpha,\lambda}(t)$ is a long memory process.   

The fractal dimension of  $B_{\alpha\beta,\lambda}(t)$ can be determined in a similar way as for  $B_{\alpha,\lambda}(t)$.
By noting that the variance of the increment process of  $B_{\alpha\beta,\lambda}(t+\tau) -B_{\alpha\beta,\lambda}(t)$
varies as $\tau^{2\alpha\beta-1}$ as $\tau \to 0$,
and using the same argument as for $B_{\alpha,\lambda}(t)$,
it can be verified that the fractal dimension of the graph of $B_{\alpha\beta,\lambda}(t)$ is a.e.\ equal to 
$\frac{5}{2} - \alpha\beta$ or $2 - H$ if $\alpha\beta = H + 1/2$.

On the other hand, the long-time asymptotic behaviour of the covariance varies as $t^{-(1+2\beta)}$,
which is independent of $\alpha$.
In contrast to FBM and TFBM $B_{\alpha,\lambda}(t)$
it is possible to separately characterize the short-time property such as fractal dimension, 
and the long-time behaviour like long-range dependence of $B_{\alpha\beta,\lambda}(t)$ by using two different parameters.

\section{Tempered Fractional Brownian Motion with Variable Tempering Parameter}
\label{sec:varpara}
\noindent
From the point of view of applications of TFBM, 
it will be useful to allow the tempering parameter to vary with time (or position).  
This section deals with TFBM with time-dependent tempering parameter $\lambda(t)$.
There exist several candidates, 
only one such generalization which appears to be most promising is considered here in details, 
and the rest are briefly discussed in the Appendix.

Recall that for TFBM, an exponential tempering with constant damping parameter  is added to the moving average representation of FBM. 
Since in applications of modelling real data, 
the tempering parameter  often depends on the sample size, time intervals, 
or the particle concentration of the heterogenous medium, etc.\ 
\cite{KullbergCastillo-Negrete2012,ZhangBaeumerReeves2010,MeerschaertSabzikarM.PhanikumarEtAl2014}.
Thus, instead of constant tempering, 
TFBM with a time (or position) dependent tempering parameter will provide a more flexible and realistic model. 
In some diffusive transport of some complex systems there exists a crossover from anomalous to normal diffusion 
which is modelled by two separate power-laws. 
Tempered anomalous diffusion model \cite{KullbergCastillo-Negrete2012,Molina-GarciaSandevSafdariEtAl2018}
with the incorporation of variable tempering parameter can provide a better description of systems with multi-scale heterogeneities.

Consider the generalization of the fractional Langevin equation 
(\ref{eq:RfOU_0010})
to variable tempering parameter:
\begin{align}
  \bigl[D+\lambda(t)\bigr]^\alpha X_{\alpha,\lambda(t)}(t) & = \eta(t) .
\label{eq:varpara_0010}
\end{align}
The case with $\alpha = 1$ and random fluctuating damping $\lambda(t)$  has been studied \cite{EabLim2018}.
The fractional case is more complicated. 
For simplicity, assume $\lambda(t)$ to be a bounded deterministic function of $t$.
Note that the operator identity (\ref{eq:RfOU_0030}) does not hold for time-dependent tempering, 
that is
 $\bigl({_aD_t} + \lambda(t)\bigr)^\alpha \neq e^{-\lambda(t)t} {_aD_t^{\alpha}} e^{\lambda(t)t}$.
Furthermore, $\bigl({_aD_t} + \lambda(t)\bigr)^\alpha$ is not well-defined for fractional $\alpha$, 
since its binomial expansion is ambiguous (see Appendix). It is possible to regard (\ref{eq:varpara_0010}) as
a pseudo-differential equation analogous to 
the treatment for the case of generalised fractional random field with variable order
\cite{Ruiz-MedinaAnhAngulo2004}.

One can define FOU with time-dependent tempering directly without associating it to fractional Langevin equation. 
Consider the Weyl FOU given by (\ref{eq:RfOU_0040}) with variable tempering parameter $\lambda(t)$ for $\alpha > 1/2$.
\begin{align}
  X_{\alpha,\lambda(t)}(t) & = \frac{1}{\Gamma(\alpha)} \int_{-\infty}^t e^{-\lambda(t)(t-u)}(t-u)^{\alpha-1} \eta(u) du.
\label{eq:varpara_0020}
\end{align}
For $t > s > 0$, its covariance can be computed as
\begin{align}
  \bigl<X_{\alpha,\lambda(t)}(t)X_{\alpha,\lambda(s)}(s)\bigr> 
            & = \frac{1}{\Gamma(\alpha)^2}
                \int_{-\infty}^{\min(t,s)} du e^{-\lambda(t)(t-u)-\lambda(s)(s-u)} (t-u)^{\alpha-1} (s-u)^{\alpha-1} \nonumber \\
            & =  \frac{e^{-\lambda_{-}(t,s)(t-s)} (t-s)^{\alpha-1/2}}{\Gamma(\alpha)\sqrt{\pi}\bigl(2\lambda_{+}(t,s)\bigr)^{\alpha-1/2}}
                K_{\alpha-1/2}\bigl(\lambda_{+}(t,s)(t-s)\bigr) .
\label{eq:varpara_0030}
\end{align}
where 
$\lambda_{+}(t,s) = \frac{\lambda(t)+\lambda(s)}{2}$
and
$\lambda_{-}(t,s) = \frac{\lambda(t)-\lambda(s)}{2}$.
The variance is
\begin{align}
    \bigl<X_{\alpha,\lambda(t)}(t)^2\bigr>  & = \frac{\Gamma(2\alpha-1)}{\bigl(\Gamma(\alpha)\bigr)^2\bigl(2\lambda(t)\bigr)^{2\alpha-1}}.
\label{eq:varpara_0040}
\end{align}
Note that when $\lambda(t)= \lambda$, 
and (\ref{eq:varpara_0030}), (\ref{eq:varpara_0040}) 
reduces respectively to (\ref{eq:RfOU_0050}), (\ref{eq:RfOU_0060}),
the covariance and variance of Weyl FOU.

For $t < s$, 
the results can be obtained by simply interchanging $t$ and $s$. 
Let
\begin{align}
  \lambda_{-}^*(t,s) & =
                       \begin{cases}
                         \lambda_{-}(t,s) & \text{for} \ t > s \\
                         \lambda_{-}(s,t) & \text{for} \ t < s 
                       \end{cases},
\label{eq:varpara_0041}
\end{align}
then for $t < s$,
\begin{align}
C_{\alpha,\lambda(t),\lambda(s)}(t,s) 
& = \frac{e^{-\lambda_{-}^*(t,s)|t-s|}}{\sqrt{\pi}\Gamma(\alpha)}
                                     \left(\frac{|t-s|}{2\lambda_{+}(t,s)}\right)^{\alpha-1/2} 
                                      K_{\alpha-1/2}\bigl(\lambda_{+}(t,s)|t-s|\bigr) .
\label{eq:varpara_0042}
\end{align}

The reduced process $B_{\alpha,\lambda(t)}(t) = X_{\alpha,\lambda(t)}(t) - X_{\alpha,\lambda(t)}(0)$
has the following covariance
\begin{gather}
  \bigl<B_{\alpha,\lambda(t)}(t)B_{\alpha,\lambda(s)}(s)\bigr> 
   = \frac{e^{-\lambda_{-}(t,s)|t-s|} |t-s|^{\alpha-1/2}}{\sqrt{\pi}\Gamma(\alpha)\bigl(2\lambda_{+}(t,s)\bigr)^{\alpha-1/2}}
                K_{\alpha-1/2}\bigl(\lambda_{+}(t,s)(t-s)\bigr) 
     + \frac{\Gamma(2\alpha-1)\bigl(2\lambda_{+}(t,s)\bigr)^{1-2\alpha}}{\Gamma(\alpha)^2} \nonumber \\
    - \frac{e^{-\lambda_{-}(t,s)t} t^{\alpha-1/2}}{\Gamma(\alpha)\sqrt{\pi}\bigl(2\lambda_{+}(t,s)\bigr)^{\alpha-1/2}}
                K_{\alpha-1/2}\bigl(\lambda_{+}(t,s)t\bigr) 
    - \frac{e^{-\lambda_{-}(t,s)s} s^{\alpha-1/2}}{\Gamma(\alpha)\sqrt{\pi}\bigl(2\lambda_{+}(t,s)\bigr)^{\alpha-1/2}}
                K_{\alpha-1/2}\bigl(\lambda_{+}(t,s)s\bigr) .
\label{eq:varpara_0050}
\end{gather}
(\ref{eq:varpara_0050}) can be expressed in the form (\ref{eq:RfOU_0140}) analogous to TFBM $B_{\alpha,\lambda}(t)$:
\begin{align}
    \bigl<B_{\alpha,\lambda(t)}(t)B_{\alpha,\lambda(s)}(s)\bigr> & = \frac{1}{2}\Bigl[c_t|t|^{2H} + c_s|s|^{2H} - c_{t-s}|t-s|^{2H}\Bigr],
\label{eq:varpara_0060}
\end{align}
with the difference that the coefficients $c_t(t,s)$, $c_s(t,s)$ and $c_{t-s}(t,s)$ depend on both $s$ and $t$,
\begin{align}
  c_t(t,s) & = \frac{2\Gamma(2\alpha-1)}{\Gamma(\alpha)^2\bigl(2\lambda_{+}(t,s)|t|\bigr)^{2\alpha-1}}
             - \frac{e^{-\lambda^*_{-}(t,s)|t|}}{\Gamma(\alpha)\sqrt{\pi}\bigl(2\lambda_{+}(t,s)|t|\bigr)^{\alpha-1/2}}
                K_{\alpha-1/2}\bigl(\lambda_{+}(t,s)|t|\bigr) ,
\label{eq:varpara_0070}
\end{align}
and $H=\alpha - \frac{1}{2}$.

In contrast TFBM $B_{\alpha,\lambda}(t)$, 
the presence of the terms $\lambda_{+}(t,s)$ and $\lambda_{-}(t,s)$ in its covariance results in non-stationarity in the increment process of $B_{\alpha,\lambda(t)}(t)$, 
and it does not satisfy the scaling property (\ref{eq:RfOU_Scaling_0010}).
However, it can be shown that both $B_{\alpha,\lambda(t)}(t)$ and $B_{\alpha,\lambda}(t)$ have similar short-time and long-time properties.

Denote by $C_{\alpha,\lambda(t),\lambda(s)}(t,s)$ the covariance (\ref{eq:varpara_0030}) of $X_{\alpha,\lambda(t)}(t)$.
Using \cite{AbramowitzStegun64}
\begin{align}
  K_\nu(z) & = 2^{\nu-1} \Gamma(\nu)z^{-\nu} + 2^{-\nu-1}\Gamma(-\nu) z^\nu 
             + \frac{z^{\nu-1}}{1-\nu}  \Gamma(\nu)z^{2-\nu} + \cdots &
         z & \to 0,
\label{eq:varpara_0080}
\end{align}
one gets
\begin{align}
  C_{\alpha,\lambda(t+\tau),\lambda(t)}(t+\tau,t) & \sim e^{-\lambda_{-}(+\tau,t)} 
                                                    \Bigl[
                                                    A \bigl(\lambda_{+}(t+\tau,t)\bigr)^{1-2\alpha}
                                                    + B \tau^{2\alpha-1}
                                                    + o\bigl(\tau^2\bigr)
                                                    \Bigr],
\label{eq:varpara_0090}
\end{align}
where $A$ and $B$ are constants.
For $\tau \to 0$, $e^{-\lambda_{-}(t+\tau,t)}  \sim 1$, 
the first term is finite and dependent on $\lambda_{+}(t,s)$ or $\lambda(t)$ and $\lambda(s)$.
The second term $\sim \tau^{2\alpha-1}$ is similar to that for TFBM with constant tempering parameter.  

For the long-time behavior, using 
\begin{align}
  K_{\nu}(z) & = \Bigl(\frac{\pi}{2}\Bigr)^{1/2} \frac{e^{-z}}{\sqrt{z}}\bigl[1+O(1/z)\bigr], &
          z & \to \infty,
\label{eq:varpara_0100}
\end{align}
one gets exponential decay for the covariance functions of both $X_{\alpha,\lambda}(t)$ and $X_{\alpha,\lambda(t)}(t)$.
Apply the same argument to $B_{\alpha,\lambda}(t)$ and $B_{\alpha,\lambda(t)}(t)$, 
gives similar short- and long-time behaviour for the covariance functions 
for these two processes.

Since the short-time properties of for both $B_{\alpha,\lambda(t)}(t)$ and $B_{\alpha,\lambda}(t)$  are similar, 
both these processes have the same fractal dimension, which is a local property.
In addition,  $B_{\alpha,\lambda(t)}(t)$ satisfies the locally self-similarity property just like $B_{\alpha,\lambda}(t)$ \cite{LimEab2019}.

Since the two TFBM processes $B_{\alpha,\lambda(t)}(t)$ and $B_{\alpha,\lambda}(t)$ have similar long time behavior, 
$B_{\alpha,\lambda(t)}(t)$ is a long memory process, just like $B_{\alpha,\lambda}(t)$.
This is expected since any finite variation in the tempering strength does not alter the overall memory characteristic of the process.

\section{Tempered Multifractional Brownian Motion}
\label{sec:TMBM}
\noindent
TFBM provides a useful model for describing geophysical flows phenomena such as wind speed modeling 
\cite{SabzikarMeerschaertChen2015,MeerschaertSabzikarM.PhanikumarEtAl2014} and diffusive transport processes 
\cite{KullbergCastillo-Negrete2012,Molina-GarciaSandevSafdariEtAl2018}.
The main attractiveness of this model is its simplicity with its properties described by a single index. 
However, situation in the real world is more complicated to be characterized by a single parameter. 
For example, in heterogenous medium with scaling that may vary with time or position, 
and the system may have variable time-dependent memory. 
Therefore, a more realistic model requires a time-dependent or position-dependent scaling exponent. 
This can be achieved by generalising TFBM to a variable index or tempered multifractional Brownian motion (TMBM). 
The extension of TMBM to include variable tempering parameter will also be discussed.

Consider first the extension of FOU to variable index or multifractional Ornsttein-Uhlenbeck process (MOU). 
TMBM can then be defined as the reduced process of MOU. 
Analogous to the generalisation of FBM to multifractional Brownian motion \cite{Levy-VehelPeltier1996,BenassiJaffardRoux1997}
one gets MOU of Weyl type by replacing the constant index $\alpha$ of Weyl FOU in (\ref{eq:RfOU_0040}) 
by a variable index $\alpha(t)$:
\begin{align}
  X_{\alpha(t),\lambda}(t) & = \frac{1}{\Gamma\bigl(\alpha(t)\bigr)} \int_{-\infty}^t e^{\lambda(t-u)} (t -u)^{\alpha(t)-1} \eta(u) du.
\label{eq:TMBM_0010} 
\end{align}
It is assumed that $\alpha(t) > 1/2$, and the H\"{o}lder continuous with $\bigl|\alpha(t) - \alpha(s)\bigr| \leq K|t-s|^\kappa$,
$K > 0$, $\kappa > 0$.
\begin{subequations}
\label{eq:TMBM_0020}
For $s <t$. the covariance of $X_{\alpha(t),\lambda}(t)$ is given by
\begin{align}
  \Bigl<X_{\alpha(t),\lambda}(t)X_{\alpha(s),\lambda}(s)\Bigr> & = \frac{e^{-\lambda(t+s)}}{\Gamma\bigl(\alpha(t)\bigr)\Gamma\bigl(\alpha(s)\bigr)} \int_{-\infty}^{\min(t,s)} (t -u)^{\alpha(t) -1} (s -u)^{\alpha(s) -1} e^{2\lambda{u}} du \nonumber \\
& = \frac{e^{-\lambda(t-s)}}{\Gamma\bigl(\alpha(t)\bigr)\Gamma\bigl(\alpha(s)\bigr)} \int_{-\infty}^s (u)^{\alpha(s) -1} (u+t- s)^{\alpha(t) -1} e^{-2\lambda{u}} du \nonumber\\
  & =  \frac{e^{-\lambda(t-s)}(t - s)^{2\alpha_{+}(t,s)-1}}{\Gamma\bigl(\alpha(t)\bigr)}
       \Psi\bigl(\alpha(s),2\alpha_{+}(t,s), 2\lambda(t-s)\bigr) ,
\label{eq:TMBM_0020a}
\end{align}
where  $\Psi(\alpha,\gamma;z)$ is the confluent hypergeometric function, which is also known as Kummer function 
(3.383 of \cite{GradshteynRyzhik2000}) and $\alpha_{+}(t,s) = \bigl(\alpha(t)+\alpha(s)\bigr)/2$.

For $s > t$, the covariance of $X_{\alpha(t),\lambda}(t)$ is given by
\begin{align}
    \Bigl<X_{\alpha(t),\lambda}(t)X_{\alpha(s),\lambda}(s)\Bigr>  & =  \frac{e^{-\lambda(s-t)}(s - t)^{2\alpha_{+}(t,s)-1}}{\Gamma\bigl(\alpha(s)\bigr)}
       \Psi\bigl(\alpha(t),2\alpha_{+}(t,s), 2\lambda(s-t)\bigr).
\label{eq:TMBM_0020b}
\end{align}  
\end{subequations}
The two cases \mbox{(\ref{eq:TMBM_0020a})} and \mbox{(\ref{eq:TMBM_0020b})}  can be combined into a single expression
\begin{align}
    \Bigl<X_{\alpha(t),\lambda}(t)X_{\alpha(s),\lambda}(s)\Bigr>  & =  \frac{e^{-\lambda|t-s|}|t - s|^{2\alpha_{+}(t,s)-1}}{\Gamma\bigl(\alpha(t\vee{s})\bigr)}
       \Psi\bigl(\alpha(t\wedge{s}),2\alpha_{+}(t,s), 2\lambda|t-s|\bigr), 
\label{eq:TMBM_0020x} 
\end{align}
where $f(t\vee{s})  = f\bigl(\max(t,s)\bigr)$ and $f(t\wedge{s})  = f\bigl(\min(t,s)\bigr)$.

There exists another possible way of defining MOU is based on the spectral representation (\ref{eq:RfOU_0080}) of FOU  
\begin{align}
  X_{\alpha(t),\lambda}(t) & = \frac{1}{\sqrt{2\pi}}\int_{-\infty}^\infty \frac{e^{-ikt}\widetilde{\eta}(k)dk}{\bigl(-ik + \lambda\bigr)^{\alpha(t)}} .
\label{eq:TMBM_0030}
\end{align}
The covariance is given by
\begin{align}
  \Bigl<X_{\alpha(t),\lambda}(t)X_{\alpha(s),\lambda}(s)\Bigr> & = \frac{1}{2\pi}\int_{-\infty}^\infty \frac{e^{-ik(t-s)}dk}{\bigl(-ik + \lambda\bigr)^{\alpha(t)}\bigl(ik + \lambda\bigr)^{\alpha(s)}}  \nonumber\\
    & = \frac{(t-s)^{\alpha_{+}(t,s)-1}}{\Gamma\bigl(\alpha(t)\bigr)(2\lambda)^{\alpha_{+}(t,s)}}
        W_{\alpha_{-}(t,s), 1/2 - \alpha_{+}(t,s) }\bigl(2\lambda(t-s)\bigr) ,
\label{eq:TMBM_0040}
\end{align}
where $W(\cdot)$ is the Whittaker function, and $\alpha_{-}(t,s) = \bigl(\alpha(t) - \alpha(s)\bigr)/2$
(see 3.384.9 of  \cite{GradshteynRyzhik2000}).
It can be shown that the two representations (\ref{eq:TMBM_0020a}) and (\ref{eq:TMBM_0040}) of MOU are equivalent \cite{LimEab2019}.

Note that just like the case of FOU, the $\lambda \to 0$ limit of the covariance and variance of 
the Weyl multifractional Ornsttein-Uhlenbeck process (MOU) diverges. 
However, the $\lambda \to 0$ limit of the reduced process of Weyl MOU is multifractional Brownian motion \cite{LimTeo2007}.

In order to obtain a TMBM which has a close semblance of TFBM, 
especially its covariance, one uses MOU of Riesz type. 
Recall that in Section 3, Riesz type FOU with two indices is introduced. 
By letting $\beta = 1$ and replacing $\alpha$ by $\alpha(t)$ in (\ref{eq:TfOU2ind_0070}) results in 
\begin{align}
  Y_{\alpha(t),\lambda}(t) & = \frac{1}{\sqrt{2\pi}}\int_{-\infty}^\infty \frac{e^{ikt}\widetilde{\eta}(k)dk}{\bigl(|k|^2 + \lambda^2\bigr)^{\alpha(t)/2}} ,
\label{eq:TMBM_0050}
\end{align}
which is the Riesz MOU. Its covariance and variance can be calculated for  $\alpha(t) > 1/2$:
\begin{align}
  \Bigl<Y_{\alpha(t),\lambda}(t)Y_{\alpha(s),\lambda}(s)\Bigr> & = \frac{1}{\sqrt{\pi}\Gamma\bigl(\alpha_{+}(t,s)\bigr)}
                                                           \left(\frac{|t-s|}{2\lambda}\right)^{\alpha_{+}(t,s) - 1/2}
                                                           K_{\alpha_{+}(t,s) - 1/2}\bigl(\lambda|t - s|\bigr) ,     
\label{eq:TMBM_0060}
\end{align}
and 
\begin{align}
  \Bigl<\bigl(Y_{\alpha(t),\lambda}(t)\bigr)^2\Bigr>  & = \frac{\Gamma\bigl(2\alpha(t)-1\bigr)}{\Bigl(\Gamma\bigl(\alpha(t)\bigr)\Bigr)^2 \bigl(2\lambda\bigr)^{2\alpha(t)-1}} .
\label{eq:TMBM_0070}
\end{align}
Note that (\ref{eq:TMBM_0050}) is a special case of the multifractional Riesz-Bessel process \cite{LimTeo2008a}.

One can define TMBM as the reduced process of Reisz MOU:                
\begin{align}
  B_{\alpha(t),\lambda}(t) & = Y_{\alpha(t),\lambda}(t) - Y_{\alpha(t),\lambda}(0),
\label{eq:TMBM_0071}
\end{align}
It is a centred Gaussian process with the following covariance and variance: 
\begin{align}
  \wideparen{C}_{\alpha(t),\lambda}(t - s) & =  \Bigl<B_{\alpha(t),\lambda}(t)B_{\alpha(s),\lambda}(s)\Bigr> \nonumber \\
                                       & = \frac{1}{\sqrt{\pi}\alpha_{+}(t,s)}
                                            \Biggl[
                                            \left(\frac{|t-s|}{2\lambda}\right)^{\alpha_{+}(t,s) - 1/2}
                                            K_{\alpha_{+}(t,s) - 1/2}\bigl(\lambda|t - s|\bigr) \nonumber \\
                                      & \qquad\qquad   - \left(\frac{|t|}{2\lambda}\right)^{\alpha_{+}(t,s) - 1/2}
                                           K_{\alpha_{+}(t,s) - 1/2}\bigl(\lambda|t|\bigr)   \nonumber \\
                                      & \qquad\qquad\qquad 
                                        - \left(\frac{|s|}{2\lambda}\right)^{\alpha_{+}(t,s) - 1/2}
                                           K_{\alpha_{+}(t,s) - 1/2}\bigl(\lambda|s|\bigr)  
                                           \Biggr]  \nonumber \\
                                      & \qquad\qquad\qquad\qquad 
                                        + \frac{\Gamma(2\alpha_{+}(t,s) - 1)}{\bigl(\Gamma(\alpha_{+}(t,s))\bigr)^2(2\lambda)^{2\alpha_{+}(t,s) - 1}}.
\label{eq:TMBM_0080}
\end{align}
and
\begin{align}
  \wideparen{\sigma}_{\alpha(t),\lambda}(t) & =  \Bigl<\bigl(B_{\alpha(t),\lambda}\bigr)^2\Bigl> \nonumber \\
                                                        & = \frac{2}{\sqrt{\pi}\Gamma\bigl(\alpha(t)\bigr)} 
                                                          \left[
                                                          \frac{\Gamma\bigl(2\alpha(t) -1\bigr)}
                                                               {\Bigl(\Gamma\bigl(\alpha(t)\bigr)\Bigr)^2(2\lambda)^{2\alpha(t)-1}}
                                                          - \left(\frac{|t|}{2\lambda}\right)^{\alpha(t) - 1/2}
                                                          K_{\alpha(t) - 1/2}\bigl(\lambda|t|\bigr)
                                                          \right],
\label{eq:TMBM_0090}
\end{align}
using $\alpha_{+}(t,s) - 1/2 = H_{+}(t,s) = \bigl(H(s)+H(t)\bigr)/2$, 
the covariance of TMFM can be expressed in the form similar to that for TFBM:
\label{eq:TMfBM_0090}
\begin{align}
  \wideparen{C}_{\alpha_{+}(t,s),\lambda}(t - s) & = \frac{1}{2} \Bigl[
                                          c_t\bigl(H_{+}(t,s)\bigr)|t|^{2H_{+}(t,s)}
                                         + c_s\bigl(H_{+}(t,s)\bigr)|s|^{2H_{+}(t,s)} \nonumber \\
                                    &  \qquad\qquad\qquad\qquad\qquad
                                         - c_{t-s}\bigl(H_{+}(t,s)\bigr)|t-s|^{2H_{+}(t,s)}
                                         \Bigr],
\label{eq:TMBM_0100}
\end{align}
with
\begin{align}
  c_t\bigl(H_{+}(t,s)\bigr) & = \frac{2\Gamma\bigl(2H_{+}(t,s)\bigr)}{\Gamma\bigl(2H_{+}(t,s) + 1/2\bigr)(2\lambda|t|)^{2H_{+}(t,s)}}
                             \nonumber \\
                           & \qquad\qquad\qquad
                             - \frac{2}{\sqrt{\pi}\Gamma\bigl(H_{+}(t,s)+1/2\bigr)}
                               \left(\frac{1}{2\lambda|t|}\right)^{H_{+}(t,s)} K_{H_{+}(t,s)}\bigl(\lambda|t|\bigr).
\label{eq:TMBM_0110}
\end{align}

Due to its variable index,  TMBM $B_{\alpha(t),\lambda}(t)$ does not satisfy the scaling property (\ref{eq:RfOU_Scaling_0010}).
Similarly, the increment process of  $B_{\alpha(t),\lambda}(t)$ is non-stationary. 
Despite the difference in the global properties between TFBM and TMBM, their short- and long-time properties are similar.

TMBM satisfies the locally asymptotically self-similar property. 
Assume $\alpha(t)$ is H\"{o}lder continuous with $|\alpha(t) - \alpha(s)| \leq K|t-s|^\kappa$ ,
and $\frac{1}{2} < \alpha(t) < \kappa + \frac{1}{2}$ for all $t$.

Consider the increment 
\begin{align}
  B_{\alpha(t),\lambda}(t) - B_{\alpha(s),\lambda}(s) & = \bigl(Y_{\alpha(t),\lambda}(t) - Y_{\alpha(s),\lambda}(s)\bigr)  
                                                - \bigl(Y_{\alpha(t),\lambda}(0) - Y_{\alpha(s),\lambda}(0)\bigr) .         
\label{eq:TMBM_0120}
\end{align}
Using a similar argument as for multifractional Riesz-Bessel process 
\cite{LimTeo2008a},
one has for $|t-s| \to 0$                         
\begin{align}
  \bigl(Y_{\alpha(t),\lambda}(0) - Y_{\alpha(s),\lambda}(0)\bigr) & = O\bigl(|t-s|^{2\kappa}\bigr) ,          
\label{eq:TMBM_0130}
\end{align}
and
\begin{align}
    \Bigl<\bigl(Y_{\alpha(t),\lambda}(t) - Y_{\alpha(s),\lambda}(s)\bigr)^2\Bigr> & = \frac{\Gamma\Bigl(\frac{1}{2} - \alpha(t)\Bigr)}
                                                                                    {2^{2\alpha(t) -1}\sqrt{\pi}\Gamma\bigl(\alpha(t)\bigr)}
                                                                                    |s - t|^{2\alpha(t)-1}
                                                                             + O\bigl(|s - t|^{2\alpha(t) - 1/2}\bigr) .     
\label{eq:TMBM_0140}
\end{align}

Consider
\begin{align}
  C_{\alpha(t)}(\epsilon;u,v) & = \Biggl<
                                \left(\frac{Y_{\alpha(t),\lambda}(t+\epsilon{u}) - Y_{\alpha(t),\lambda}(t)}{\epsilon^{\alpha(t) - 1/2}}\right)
                                \left(\frac{Y_{\alpha(t),\lambda}(t+\epsilon{v}) - Y_{\alpha(t),\lambda}(t)}{\epsilon^{\alpha(t) - 1/2}}\right)
                                \Biggr> \nonumber \\
     & = \frac{1}{2\epsilon^{2\alpha(t)-1}}
          \biggl[
              \Bigl<\bigl(Y_{\alpha(t),\lambda}(t+\epsilon{u}) - Y_{\alpha(t),\lambda}(t)\bigr)^2\Bigr> \nonumber \\
       & \qquad\qquad\qquad
         + \Bigl<\bigl(Y_{\alpha(t),\lambda}(t+\epsilon{v}) - Y_{\alpha(t),\lambda}(t)\bigr)^2\Bigr>  \nonumber \\
       & \qquad\qquad\qquad\qquad
         - \Bigl<\bigl(Y_{\alpha(t),\lambda}(t+\epsilon{v}) - Y_{\alpha(t),\lambda}(t+\epsilon{v})\bigr)^2\Bigr>
          \biggr].
\label{eq:TMBM_0150}
\end{align}
For $\epsilon \to 0$ one has from (\ref{eq:TMBM_0140})
\begin{align}
  \Bigl<\bigl(B_{\alpha(t),\lambda}(t+\epsilon{u}) - B_{\alpha(t),\lambda}(t))\bigr)^2\bigr> 
         & = \frac{\Gamma\bigl(1/2 - \alpha(t)\bigr)}{2^{2\alpha(t)-1}\sqrt{\pi}\Gamma\bigl(\alpha(t)\bigr)}
           \bigl(\epsilon|u|\bigl)^{2\alpha(t)-1} + O\bigl(\epsilon^{\kappa+ \alpha(t)-1/2}\bigr).
\label{eq:TMBM_0160}
\end{align}
Therefore,
\begin{align}
  \lim_{\epsilon\to 0} C_{\alpha(t)}(\epsilon;u,v) & = \lim_{\epsilon\to 0} \frac{1}{2\epsilon^{2\alpha(t)-1}}
                                                   \Biggl[
                                                      \frac{\Gamma\bigl(1/2 - \alpha(t)\bigr)}{2^{2\alpha(t)-1}\sqrt{\pi}\Gamma\bigl(\alpha(t)\bigr)}
                                                     \Bigl(\bigl(\epsilon|u|\bigr)^{2\alpha(t) -1} + \bigl(\epsilon|v|\bigr)^{2\alpha(t) -1}\Bigr)
                                                 \nonumber \\
           & - \frac{\Gamma\bigl(1/2 - \alpha(t+\epsilon{v})\bigr)}{2^{2\alpha(t)-1}\sqrt{\pi}\Gamma\bigl(\alpha(t+\epsilon{v})\bigr)}
             \bigl(\epsilon|u-v|\bigr)^{2\alpha(t+\epsilon{v}) -1} + O\bigl(\epsilon^{\kappa+\min\bigl(\alpha(t),\alpha(t+\epsilon{v})\bigr)-1/2}\bigr)
             \Biggr]  \nonumber \\
     & = \frac{\Gamma\bigl(1/2 - \alpha(t)\bigr)}{2^{2\alpha(t)-1}\sqrt{\pi}\Gamma\bigl(\alpha(t)\bigr)}
       \Bigl(|u|^{2\alpha(t) -1} + |v|^{2\alpha(t) -1} - |u-v|^{2\alpha(t) -1}\Bigr), 
\label{eq:TMBM_0170}
\end{align}
where H\"{o}lder continuity of $\alpha(t)$ implies
\begin{align}
  \frac{\Gamma\bigl(1/2 - \alpha(t+\epsilon{v})\bigr)}{2^{2\alpha(t)-1}\sqrt{\pi}\Gamma\bigl(\alpha(t+\epsilon{v})\bigr)}
  \bigl(\epsilon|u-v|\bigr)^{2\alpha(t+\epsilon{v}) -1}  
  &  \rightarrow  \frac{\Gamma\bigl(1/2 - \alpha(t)\bigr)}{2^{2\alpha(t)-1}\sqrt{\pi}\Gamma\bigl(\alpha(t)\bigr)} 
    \bigl(\epsilon|u-v|\bigr)^{2\alpha(t) -1}  .
\label{eq:TMBM_0180}
\end{align}
Note that (\ref{eq:TMBM_0170}) is the covariance of multifractional Brownian motion. 
Thus, TMBM is locally asymptotically self-similar, 
and its tangent process at a point $t_\circ$ is the multifractional Brownian motion with Hurst index 
$H(t_\circ) = \alpha(t_\circ) - 1/2$.

Another local property is the fractal dimension of the graph of TMBM. 
With probability one, 
the Hausdorff dimension of the graph of the TMBM $Y_{\alpha(t),\lambda}(t)$ indexed by $\alpha(t)$ over the intercal
$I \in \mathbb{R}$ is $\frac{5}{2} - m_I\bigl[\alpha(t)\bigr]$, where $m_I\bigl[\alpha(t)\bigr] = \min\{\alpha(t); t \in I\}$.
First note that the leading local behaviours of the variance of increments of multifractional Brownian motion and TMBM are the same. 
Since only this leading behavior is used to arrive at the local properties of the TMBM, 
one can easily infer that all these properties of the TMBM holds verbatim for the MBM if we identify 
$H(t)$ with $\alpha(t) - 1/2$ \cite{BenassiCohenIstas2003}.

Finally, one would like to consider TMBM with a time-dependent tempered parameter $\lambda(t)$.
Such a process $B_{\alpha(t),\lambda(t)}(t)$ can provide a more flexible model 
for systems which require both the variable scaling index and tempered parameter, 
for example turbulence flow in a heterogenous medium.

In previous section, TFBM with variable tempering $B_{\alpha,\lambda(t)}(t)$ is formulated in terms of reduced process of Weyl type FOU 
$X_{\alpha,\lambda(t)}(t)$.
It is also known that $Y_{\alpha,\lambda(t)}(t)$ based on FOU of Riesz type does not lead to a viable TFBM.
Therefore, one would expect the reduced process of Weyl MOU given by 
(\ref{eq:TMBM_0010})
as a suitable candidate for TMBM with variable tempering parameter.

Let $X_{\alpha(t),\lambda(t)}(t)$ be the Weyl MOU with variable tempering parameter defined by 
\begin{align}
  X_{\alpha(t),\lambda(t)}(t) & = \frac{1}{\Gamma\bigl(\alpha(t)\bigr)}\int_{-\infty}^t e^{-\lambda(t)(t-u)}(t-u)^{\alpha(t)-1} \eta(u) du.
\label{eq:TMBM_1010}
\end{align}

Its covariance can be calculated for $t > s$ and $\alpha(t) > 1/2$,
\begin{subequations}
\label{eq:TMBM_1020}
\begin{align}
  \bigl<X_{\alpha(t),\lambda(t)}(t)X_{\alpha(s),\lambda(s)}(s)\bigr> 
  & = \frac{e^{\lambda(t)(t-s)}(t-s)^{2\alpha_{+}(t,s)-1}}{\Gamma\bigl(\alpha(t)\bigr)}
       \Psi\bigl(\alpha(s),2\alpha_{+}(t,s),2\lambda_{+}(t,s)(t-s)\bigr).
\label{eq:TMBM_1020a}
\end{align}
Note that this covariance is symmetric with respect to $t$ and $s$. 
In the case $s > t > 0$,
\begin{align}
  \bigl<X_{\alpha(t),\lambda(t)}(t)X_{\alpha(s),\lambda(s)}(s)\bigr> 
  & = \frac{e^{\lambda(t)(s-t)}(s-t)^{2\alpha_{+}(t,s)-1}}{\Gamma\bigl(\alpha(s)\bigr)}
       \Psi\bigl(\alpha(t),2\alpha_{+}(t,s),2\lambda_{+}(t,s)(s-t)\bigr).
\label{eq:TMBM_1020b}
\end{align}  
\end{subequations}
These two cases can be combined into single expression
\begin{align}
  \bigl<X_{\alpha(t),\lambda(t)}(t)X_{\alpha(s),\lambda(s)}(s)\bigr> 
  & = \frac{e^{\lambda(t\vee{s})(t-s)}(t-s)^{2\alpha_{+}(t,s)-1}}{\Gamma\bigl(\alpha(t\vee{s})\bigr)}
       \Psi\bigl(\alpha(t\wedge{s}),2\alpha_{+}(t,s),2\lambda_{+}(t,s)(t-s)\bigr).
\label{eq:TMBM_1030}
\end{align}
for $t,s > 0$, and $t\wedge{s} = \min(t,s)$, $t\vee{s} = max(t,s)$.

The variance is
\begin{align}
\bigl<X_{\alpha(t),\lambda(t)}(t)X_{\alpha(t),\lambda(t)}(t)\bigr> 
         & = \frac{\Gamma\bigl(2\alpha(t)-1\bigr)}{\Gamma\bigl(\alpha(t)\bigr)^2 \bigl(2\lambda(t)\bigr)^{2\alpha(t)-1}}, &
         & \text{for} \ \alpha(t) > \frac{1}{2}.
\label{eq:TMBM_1040}
\end{align}
Just like in the case of constant tempering, the $\lambda(t) \to 0$ limits for both the covariance and variance 
diverge.

TMBM with variable tempering denoted by $B_{\alpha(t),\lambda(t)}(t)$ is defined by the reduced process of 
$X_{\alpha(t),\lambda(t)}(t)$:
\begin{align}
  B_{\alpha(t),\lambda(t)}(t) & = X_{\alpha(t),\lambda(t)}(t) - X_{\alpha(t),\lambda(t)}(0).
\label{eq:TMBM_1050}
\end{align}
This is Gaussian centred process with covariance
\begin{align}
  \bigl<B_{\alpha(t),\lambda(t)}(t)B_{\alpha(s),\lambda(s)}(s)\bigr> 
  & = \frac{e^{\lambda(t)(t-s)}(t-s)^{2\alpha_{+}(t,s)-1}}{\Gamma\bigl(\alpha(t)\bigr)}
      \Psi\bigl(\alpha(s),2\alpha_{+}(t,s),2\lambda_{+}(t,s)(t-s)\bigr)  \nonumber \\
  & \qquad - \frac{e^{\lambda(t)t}(t)^{2\alpha_{+}(t,s)-1}}{\Gamma\bigl(\alpha(t)\bigr)}
       \Psi\bigl(\alpha(s),2\alpha_{+}(t,s),2\lambda_{+}(t,s)t)\bigr) \nonumber \\
  & \qquad\qquad - \frac{e^{\lambda(s)s}(s)^{2\alpha_{+}(t,s)-1}}{\Gamma\bigl(\alpha(s)\bigr)}
       \Psi\bigl(\alpha(t),2\alpha_{+}(t,s),2\lambda_{+}(t,s)s)\bigr) \nonumber \\
  & \qquad\qquad\qquad + \frac{\Gamma\bigl(\alpha_{+}(t,s) - 1\bigr)}{\Gamma\bigl(t\bigr)\Gamma\bigl(s\bigr)\bigl(2\lambda_{+}(t,s)\bigr)^{2\alpha_{+}(t,s)-1}}.
\label{eq:TMBM_1051}
\end{align}
In order to express the covariance of the TMBM with variable tempering parameter  
$B_{\alpha(t),\lambda(t)}(t)$
in the similar form as TFBM $B_{\alpha,\lambda}(t)$,
the following notations are introduced.
\begin{subequations}
\label{eq:TMBM_1060}
\begin{align}
  f_{t,s}(u\vee{v}) & = \begin{cases}
                        f(t) & \text{if} \ u > v \\
                        f(s) & \text{if} \ u < v 
                      \end{cases} 
\label{eq:TMBM_1060a}\\
  f_{t,s}(u\wedge{v}) & = \begin{cases}
                        f(s) & \text{if} \ u > v \\
                        f(t) & \text{if} \ u < v 
                      \end{cases} 
\label{eq:TMBM_1060b}
\end{align}  
\end{subequations}
The covariance (\ref{eq:TMBM_1051}) 
expressed in the form similar to (\ref{eq:RfOU_0140}) is given by
\begin{align}
  \wideparen{C}_{t,s}(t,s) & = c_{t,s}(t,0)|t|^{2H_{+}(t,s)} + c_{t,s}(0,s)|s|^{2H_{+}(t,s)} - c_{t,s}(t,s)|t-s|^{2H_{+}(t,s)},
\label{eq:TMBM_1130}  
\end{align}
where the coefficients $c(t,s)$ are given by 
\begin{align}
  c_{t,s}(u,v) 
         & = \frac{\Gamma\bigl(2H_{+}(t,s)\bigr)}{\Gamma(\alpha(t))\Gamma(\alpha(s))|u-v|^{2H_{+}(t,s)}
             \bigl(2\lambda_{+}(t,s)\bigr)^{2H_{+}(t,s)}} \nonumber \\
         & \qquad -  \frac{e^{-\lambda_{t,s}(u\vee{v})|u-v|}}{\Gamma(\alpha_{t,s}(u\vee{v}))} 
                               \Psi\bigl(\alpha_{t,s}(u\wedge{v}),2\alpha_{+}(t,s),2\lambda_{+}(t,s)|u-v|\bigr).
\label{eq:TMBM_1140}  
\end{align}

Finally, the properties of $B_{\alpha(t),\lambda(t)}(t)$ are briefly considered.  
Just like TMBM $B_{\alpha(t),\lambda}(t)$,
the scaling and stationary increment properties do not hold for  $B_{\alpha(t),\lambda(t)}(t)$.
Recall the short-  and long-time properties of  
$B_{\alpha,\lambda(t)}(t)$ are the same as  $B_{\alpha,\lambda}(t)$;
a similar situation exists for $B_{\alpha(t),\lambda(t)}(t)$ and $B_{\alpha,\lambda}(t)$.
By comparing the covariance functions of $B_{\alpha(t),\lambda(t)}(t)$ and $B_{\alpha(t),\lambda}(t)$
term by term, one sees that the difference in the exponential terms 
$e^{\lambda(t)(t-s)}$ just like $e^{\lambda(t-s)}$, 
does not change the short- and long-time behaviour of the covariance function. 
Thus, both $B_{\alpha(t),\lambda(t)}(t)$ and $B_{\alpha(t),\lambda}(t)$ have the similar short- and long-time properties.
In particular, $B_{\alpha(t),\lambda(t)}(t)$ has the same fractal dimension as $B_{\alpha(t),\lambda}(t)$;
and it is a long-memory process.

\section{Concluding Remarks}
\label{sec:conclude}
\noindent
By considering TFBM and TMBM in terms of reduced processes of FOU and MOU provides some advantages. 
First, it allows one to use the known results of FOU and MOU to study the properties of TFBM and TMBM. 
They also facilitate various generalizations of TFBM in a more direct way. 
Note that TFBM and its various generalizations have the same local behavior as their corresponding FOU and MOU processes. 
In particular, TFBM and TMBM behave like FBM and MBM respectively in the small-time scales, 
However, their global properties are not the same. 
FOU with single index and two indices, 
and MOU are all short memory processes. 
On the other hand, their reduced processes (or TFBM and the corresponding generalizations) have long-memory dependence.

TFBM with two indices is a useful addition to the short list of biparametric processes. 
which have the nice property of separate characterization of the local property of fractal dimension, 
and the global property of long-range dependence. 
The decoupling of these two independent properties provides a useful stochastic model for many physical, geological, biological, 
and socio-economic systems. 
Similarly, TMBM permits modeling of systems with variable fractal dimension and memory. 
TFBM and TMBM with time-dependent tempering parameter provide models which are more realistic and flexible for turbulence in geophysical flows, 
transport processes in heterogenous media, etc. 
It is hoped that the generalisations discussed in this paper can lead to more robust and optimal models for various natural and man-made systems.


\pagebreak
\appendix
\renewcommand{\theequation}{\Alph{section}\arabic{equation}}

\renewcommand{\theequation}{A\arabic{equation}}
\section*{Appendix: Candidates for TFBM with Time-dependent Tempering Parameter}
\label{sec:TimeLExp}
\noindent
Consider first the fractional Langevin equation with time-dependent tempering parameter
\begin{align}
  \bigl[D+\lambda(t)\bigr]^\alpha X_{\alpha,\lambda(t)}(t) & = \eta(t).
\label{eq:TimeLExp_0010}
\end{align}
This equation is formal and not well-defined in the usaul sense, 
as its binomial expansion is not unique, as can be shown below. 
First note that
\begin{align}
  \bigl[D+\lambda(t)\bigr]^\alpha & = D^\alpha\bigl[1+I\lambda(t)\bigr]^\alpha = D^\alpha\sum_{n=0}^\infty\binom{\alpha}{n} \bigl[I\lambda(t)\bigr]^n,
\label{eq:TimeLExp_0020}
\end{align}
where
\begin{align}
  \bigl[I\lambda(t)\bigr]^n & = \int_{-\infty}^t du_1 \lambda(u_1) \int_{-\infty}^{u_1} du_2 \lambda(u_2)  \cdots  \int_{-\infty}^{u_{n-1}} du_n \lambda(u_n)f(u_n).
\label{eq:TimeLExp_0030}
\end{align}
The solution of (\ref{eq:TimeLExp_0010}) based on (\ref{eq:TimeLExp_0020}) is
\begin{align}
  X(t) & = \bigl[1+I\lambda(t)\bigr]^{-\alpha} I^\alpha\eta(t) = \sum_{n=0}^\infty\binom{-\alpha}{n} \bigl[I\lambda(t)\bigr]^n I^\alpha \eta(t).
\label{eq:TimeLExp_0040}
\end{align}
On the other hand, one has
\begin{align}
  \bigl[D+\lambda(t)\bigr]^\alpha & = \bigl[1+\lambda(t)I\bigr]^\alpha D^\alpha = \sum_{n=0}^\infty\binom{\alpha}{n} \bigl[\lambda(t)I\bigr]^n D^\alpha,
\label{eq:TimeLExp_0050}
\end{align}
with
\begin{align}
  \bigl[\lambda(t)I\bigr]^n & = \lambda(t)\int_{-\infty}^t du_1 \lambda(u_1) \int_{-\infty}^{u_1} du_2 \lambda(u_2)  \cdots  \int_{-\infty}^{u_{n-1}} du_n f(u_n).
\label{eq:TimeLExp_0060}
\end{align}

This gives solution of (\ref{eq:TimeLExp_0010}) as
\begin{align}
  X(t) & = I^\alpha\bigl[1+\lambda(t)I\bigr]^{-\alpha} \eta(t) =  I^\alpha\sum_{n=0}^\infty\binom{-\alpha}{n} \bigl[\lambda(t)I\bigr]^n \eta(t).
\label{eq:TimeLExp_0070}
\end{align}
Thus, (\ref{eq:TimeLExp_0010}) does not appear to be a suitable candidate for the extension of TFBM to variable tempering.

Following the approach in \cite{Ruiz-MedinaAnhAngulo2004},
one can consider (\ref{eq:TimeLExp_0010}) as a pseudo-differential equation:
\begin{align}
  \bigl\llbracket D_t+\lambda(t)\bigr\rrbracket^\alpha  X_{\lambda(t)}(t) & := \int_{-\infty}^\infty K_{\lambda(t)}(t-u)X_{\lambda(t)}(u) du,
\label{eq:TimeLExp_0040x}
\end{align}
where the derivative in $\llbracket \cdot \rrbracket$ denotes pseudo-differential operator, and 
\begin{align}
  K_{\lambda(t)}(t-u) & = \frac{1}{\sqrt{2\pi}}\int_{-\infty}^\infty e^{ik(t-u)} \bigl(ik + \lambda(t)\bigr)^\alpha dk.
\label{eq:TimeLExp_0050x}
\end{align}
The fundamental solution is given by 
\begin{align}
  G_{\lambda(t)}(t-u) & = \frac{1}{\sqrt{2\pi}}\int_{-\infty}^\infty \frac{e^{ik(t-u)}dk}{\bigl(ik + \lambda(t)\bigr)^\alpha} .
\label{eq:TimeLExp_0060x}
\end{align}
Thus the solution of (\ref{eq:TimeLExp_0040x}) is
\begin{align}
  X_{\alpha,\lambda(t)}(t) & = \int_{-\infty}^\infty   G_{\lambda(t)}(t-u)\eta(u( du  \nonumber \\
                       & = \frac{1}{\sqrt{2\pi}}\int_{-\infty}^\infty \frac{e^{ikt}\widetilde{\eta}(k)dk}{\bigl(ik + \lambda(t)\bigr)^\alpha} ,
\label{eq:TimeLExp_0061}
\end{align}
where $\widetilde{\eta}(k)$ is the Fourier transform of $\eta(t)$.

For $t > s$, the covariance is given by
\begin{align}
    C_{\lambda(t),\lambda(s)}(t,s)   & = \frac{1}{\Gamma(\alpha)^2}
                                         \int_{-\infty}^s du
                                         e^{-\lambda(t)(t-u)-\lambda(s)(s-u)}(t-u)^{\alpha-1}(s-u)^{\alpha-1} \nonumber \\
                                     & = \frac{e^{-\lambda(t)(t-s)}}{\Gamma(\alpha)^2}
                                            \int_0^\infty du
                                            e^{-2\lambda_{+}(t,s)u}
                                            (t-s+u)^{\alpha-1} u^{\alpha-1} .
\label{eq:TimeLExp_C0020}
\end{align}
Form \citet{GradshteynRyzhik2000}, page 348, \#3.383.8:
\begin{align}
  \int_0^\infty x^{\nu-1} (x+\beta)^{\nu-1} e^{-\mu{x}} dx & = \frac{1}{\sqrt{\pi}} \left(\frac{\beta}{\mu}\right)^{\nu-\frac{1}{2}} e^{\frac{\beta\mu}{2}}\Gamma(\nu)
                                                         K_{\frac{1}{2}-\nu}\left(\frac{\beta\mu}{2}\right)  \label{eq:TimeLExp_C0030}  \\
                                                     & \qquad\qquad\qquad\qquad \Bigl[|\arg{\beta}| < \pi, Re{\mu} > 0, Re{\nu} > 0\Bigr].    \nonumber     
\end{align}
Let
\begin{align}
  \beta & = t-s, &
  \mu   & = 2\lambda_{+}(t,s), &
  \nu   & = \alpha.
\label{eq:TimeLExp_C0040}  
\end{align}
and using the fact that $K_{1/2-\alpha}(x) = K_{\alpha-1/2}(x)$, the covariance 
(\ref{eq:TimeLExp_C0020})
becomes:
\begin{align}
  C_{\lambda(t),\lambda(s)}(t,s) & = \frac{e^{-\lambda(t)(t-s)}}{\Gamma(\alpha)^2}
                                                       \frac{1}{\sqrt{\pi}} 
                                                       \left(\frac{t-s}{2\lambda_{+}(t,s)}\right)^{\alpha-1/2} e^{\lambda_{+}(t,s)(t-s)}\Gamma(\alpha)
                                                       K_{\alpha-1/2}\bigl(\lambda_{+}(t,s)(t-s)\bigr).
\label{eq:TimeLExp_C0050}  
\end{align}
Using
\begin{align}
  -\lambda(t)(t-s) +\lambda_{+}(t,s)(t-s) & = \Bigl[-\lambda(t) + \frac{1}{2}\bigl(\lambda(t) + \lambda(s)\bigr)\Bigr](t-s)
                                           = -\frac{1}{2}\bigl(\lambda(t) - \lambda(s)\bigr)(t-s),
\label{eq:TimeLExp_C0060}  
\end{align}
and let
\begin{align}
  \lambda_{-}(t,s) & =  \frac{1}{2}\bigl(\lambda(t) - \lambda(s)\bigr),
\label{eq:TimeLExp_C0070}  
\end{align}
gives
\begin{align}
  C_{\lambda(t),\lambda(s)}(t,s) & = \frac{e^{-\lambda_{-}(t,s)(t-s)}}{\sqrt{\pi}\Gamma(\alpha)}
                                                       \left(\frac{t-s}{2\lambda_{+}(t,s)}\right)^{\alpha-1/2} 
                                                       K_{\alpha-1/2}\bigl(\lambda_{+}(t,s)(t-s)\bigr)
\label{eq:TimeLExp_C0070}  
\end{align}

For $t < s$, one gets the result by interchanging between t and s. 
Using
\begin{align}
  \lambda_{-}^*(t,s) & =
                       \begin{cases}
                         \lambda_{-}(t,s) & \text{for} \ t > s \\
                         \lambda_{-}(s,t) & \text{for} \ t < s 
                       \end{cases}
\label{eq:TimeLExp_C0080}  
\end{align}
one gets
\begin{align}
  C^{++}_{\lambda(t),\lambda(s)}(t,s) & = \frac{e^{-\lambda_{-}^*(t,s)|t-s|}}{\sqrt{\pi}\Gamma(\alpha)}
                                     \left(\frac{|t-s|}{2\lambda_{+}(t,s)}\right)^{\alpha-1/2} 
                                      K_{\alpha-1/2}\bigl(\lambda_{+}(t,s)|t-s|\bigr)
\label{eq:TimeLExp_C0090}  
\end{align}
This is just the covariance of Weyl FOU with time-dependent tempering parameter (\ref{eq:varpara_0030}).

Another possible way to define FOU process with variable tempering is 
\begin{align}
     Y_{\alpha,\lambda(t)}(t) & = \frac{1}{\sqrt{2\pi}}\int_{-\infty}^\infty \frac{e^{ikt}\widetilde{\eta}(k)dk}{\bigl[|k|^2+\lambda(t)^2\bigr]^{\alpha/2}}.
\label{eq:TimeLExp_0090}
\end{align}
Again following \cite{Ruiz-MedinaAnhAngulo2004} one can regard (\ref{eq:TimeLExp_0090}) 
as the solution to the following pseudo-differential equation
\begin{align}
    \bigl\llbracket-\mathbf{D}_t^2+\lambda(t)^2\bigr\rrbracket^{\alpha/2} Y_{\alpha,\lambda(t)}(t) & = \eta(t).
\label{eq:TimeLExp_0080a}
\end{align}
$-\mathbf{D}_t^2$ is the one-dimensional Riesz derivative (or one-dimensional Laplacian operator $-\Delta$)
defined by (\ref{eq:TfOU2ind_0060}).
(\ref{eq:TimeLExp_0080a}) can be written as
\begin{align}
      \bigl\llbracket\mathbf{D}_t+\lambda(t)\bigr\rrbracket^{\alpha/2}\bigl\llbracket-\mathbf{D}_t+\lambda(t)\bigr\rrbracket^{\alpha/2} Y_{\alpha,\lambda(t)}(t) & = \eta(t).
\label{eq:TimeLExp_0100}
\end{align}
Its solution is
\begin{align}
  Y_{\alpha,\lambda(t)}(t) & = \frac{1}{\Gamma(\alpha)^2} 
                          \int_{-\infty}^\infty dv \eta(v)
                         \int_{\max(t,v)}^\infty du e^{+\lambda(t)(t-2u+v)}  (u-t)^{\alpha/2}(u-v)^{\alpha/2} 
\end{align}
The solution can be expressed as
\begin{align}
  Y_{\alpha,\lambda(t)}(t) & = G^{+}_{\alpha/2,\lambda(t)}(t)\ast G^{-}_{\alpha/2,\lambda(t)}(t) \ast \eta(t)
                         = G^{-}_{\alpha/2,\lambda(t)}(t)\ast G^{+}_{\alpha/2,\lambda(t)}(t) \ast \eta(t) ,   
\label{eq:TimeLExp_K0050}
\end{align}
where 
\begin{align}
  \widetilde{G}^\pm_{\alpha,\lambda(t)}(k) & = \frac{1}{\bigl(\pm{ik} + \lambda(t)\bigr)^\alpha}.
\label{eq:TimeLExp_K0040}
\end{align}

The covariance of $Y_{\alpha,\lambda(t)}(t)$ is given by
\begin{align}
  C_{\lambda(t),\lambda(s)}(t,s) & = \frac{1}{2\pi} \int_{-\infty}^\infty dk\int_{-\infty}^\infty dk^\prime 
              \frac{e^{ikt_ik^\prime{s}} \bigl<\widetilde{\eta}(t)\widetilde{\eta}(s)\bigr>}
                   {\bigl(|k|^2 + \lambda(t)^2\bigr)^{\alpha/2}\bigl(|k^\prime|^2 + \lambda(s)^2\bigr)^{\alpha/2}} \nonumber \\
          & = \frac{1}{2\pi} \int_{-\infty}^\infty dk
              \frac{e^{ik|t-s|}}
                   {\bigl(|k|^2 + \lambda(t)^2\bigr)^{\alpha/2}\bigl(|k|^2 + \lambda(s)^2\bigr)^{\alpha/2}},
\label{eq:TimeLExp_C2010}  
\end{align}
which cannot be evaluated in a closed analytic form. 
Thus, the process $Y_{\alpha,\lambda(t)}(t)$  given by (\ref{eq:TimeLExp_0090}) does not lead to a simple FOU process with time-dependent tempering parameter.


\pagebreak

\bibliographystyle{unsrtnat}

\bibliography{biblioTFBM}
\end{document}